\documentclass[12pt]{article}
\usepackage[usenames,dvipsnames]{color}
\usepackage{xr,graphicx,amsthm,amsmath,amssymb,xargs,srcltx,bbm,todonotes,stmaryrd,framed,a4wide,accents,sectsty,marvosym,mathrsfs}
\usepackage[shortlabels]{enumitem}
\usepackage{framed}
\usepackage{pdfsync}
\usepackage{extarrows}

\usepackage[colorlinks]{hyperref}

\newcommandx{\conditiondhsumpaper}[2]{\ensuremath{\mathcal{S}(#1,#2)}} % condition for
\newcommandx{\conditiondhsumstrongerpaper}[2]{\ensuremath{\mathcal{S}_+(#1,#2)}} % condition fo
\newcommandx{\conditiondhsumstrongergammapaper}[3]{\ensuremath{\mathcal{S}^{(#3)}(#1,#2)}} \newcommandx{\conditiondhsumsdouble}[4]{\ensuremath{\mathcal{S}_{#3,#4}(#1,#2)}}

 \newcommand{\rom}[1]{
  \textup{\lowercase\expandafter{\romannumeral#1}} }

\newcommand{\Rom}[1]{
  \textup{\uppercase\expandafter{\romannumeral#1}} }

\subsubsectionfont{\itshape}

\usepackage{cleveref}

\newtheorem{theorem}{Theorem}[subsection]
\newtheorem{proposition}[theorem]{Proposition}
\newtheorem{lemma}[theorem]{Lemma}
\newtheorem{corollary}[theorem]{Corollary}
\newtheorem{definition}[theorem]{Definition}
\newtheorem{hypothesis}[theorem]{Assumption}
\newtheorem{example}[theorem]{Example}
\newtheorem{remark}{Remark}[section]

\crefname{hypothesis}{assumption}{assumptions}

\crefname{lemma}{Lemma}{Lemmas}
\crefname{theorem}{Theorem}{Theorems}
\crefname{proposition}{Proposition}{Propositions}
\crefname{exercise}{Exercise}{Exercises}
\crefname{corollary}{Corollary}{Corollaries}

\numberwithin{equation}{section}

\newcommandx\TEP[1][1=]{\mathbb{G}^{#1}} % Univariate TEP

\newcommandx\tedcluster[2][1=n,2=\dhinterseq]{{\widetilde{\boldsymbol\nu}}^*_{#1,#2}} % cluster empirical measure
\newcommandx\tedclustersl[2][1=n,2=\dhinterseq]{{\widetilde{\boldsymbol\mu}}^*_{#1,#2}} % cluster empirical
\newcommandx\tedclusterboot[3][1=n,2=\dhinterseq,3=\xi]{{\widetilde{\boldsymbol\nu}}^*_{#1,#2,#3}}
\newcommandx\tedclusterslboot[3][1=n,2=\dhinterseq,3=\xi]{{\widetilde{\boldsymbol\mu}}^*_{#1,#2,#3}} \newcommandx\tedclusterefron[3][1=n,2=\dhinterseq,3={\rm boot}]{{\widetilde{\boldsymbol\nu}}^*_{#1,#2,#3}}
\newcommandx\tedclusterslefron[3][1=n,2=\dhinterseq,3={\rm boot}]{{\widetilde{\boldsymbol\mu}}^*_{#1,#2,#3}}
\newcommandx\tedclusterindep[2][1=n,2=\dhinterseq]{{\widetilde{\boldsymbol\nu}}^\indep_{#1,#2}} % TED for cluster functionals block estimator, independent blocks
\newcommandx\tedclustermdep[2][1=n,2=\dhinterseq]{{\widetilde{\boldsymbol\nu}}^{*(m)}_{#1,#2}} % TED for cluster functionals block estimator m-dep approximation - using \tedcluster with ${\cdot}^{(m)}$ creates a problem
\newcommandx\tedclusterrandom[2][1=n,2=\dhinterseq]{{\widehat{\boldsymbol\nu}}^*_{#1,#2}}

\newcommand\constant{\mathrm{cst}}

\newcommandx\envelope[1][1={\mathbf H}]{{\mathbf {#1}}}

\newcommand\canditheta{{\vartheta}}
\newcommandx\stoploss[1][1={\rm stoploss}]{\theta_{#1}}
\newcommandx\stoplossest[1][1={{\rm stoploss},n}]{\widetilde\theta_{#1}}
\newcommandx\stoplossestk[1][1={{\rm stoploss},n,\statinterseq}]{\widehat\theta_{#1}}
\newcommandx\largedeviation[1][1={\rm largedev}]{\theta_{#1}}
\newcommandx\largedeviationest[1][1={{\rm largedev},n}]{\widetilde\theta_{#1}}
\newcommandx\largedeviationestk[1][1={{\rm largedev},n,\statinterseq}]{\widehat\theta_{#1}}
\newcommandx\ruin[1][1={\rm ruin}]{\theta_{#1}}
\newcommandx\ruinest[1][1={{\rm ruin},n}]{\widetilde\theta_{#1}}
\newcommandx\ruinestk[1][1={{\rm ruin},n,\statinterseq}]{\widehat\theta_{#1}}

\newcommandx\clfunc[1][1=h]{#1} % notation for cluster functionals
\newcommandx\anticl[2][2=\dhinterseq]{\conditionS(#1,#2)}

\newcommandx\anticlpsi[3][1=\tepseq,2=\dhinterseq,3=\psi]{\conditionS(#1,#2,#3)}

\newcommandx{\norm}[2][1=]{\left|#2\right|_{#1}} %%%% arbitrary norm on a euclidean space
\newcommandx{\lpnorm}[3][1=,3=]{\|#2\|_{#1}^{#3}}

\newcommandx{\matrixnorm}[2][2=]{\left\|#1\right\|_{#2}} %%%% matrix norm
\newcommandx{\matrixnormseries}[3][2=,3=]{\left\|#1\right\|_{#2,#3}} %%%% matrix norm
\newcommandx{\linftynorm}[2][1=]{\|#2\|_{#1}} %%%% supnorm in $\ellinfty$ of $A$
\newcommandx{\anynorm}[2][2=]{\left|#1\right|_{#2}} %%%% supnorm
\newcommandx{\supnormclass}[3][3=]{\left\|#1\right\|_{#2}^{#3}} %%%% supnorm on a function class (for entropy)

\newcommandx{\sphere}[2][1=]{\mathbb{S}_{#1}^{#2}}
\newcommandx{\oball}[3][3=]{B_{#3}(#1,#2)} %%%%% open ball with center #1 and radius #2  ; #3 is the norm
\newcommandx{\cball}[3][3=]{\overline{B}_{#3}(#1,#2)} %%%%% closed ball with center #1 and radius #2  ; #3 is the norm
\newcommandx{\coball}[3][3=]{B_{#3}^c(#1,#2)} %%%%% complement of the open ball with center #1 and radius #2
\newcommandx{\ccball}[3][3=]{\overline{B}_{#3}^c(#1,#2)} %%%%% complement of the closed ball with center #1 and radius #2
\newcommandx\cone[1][1=C]{\mathcal{#1}}

\newcommandx\conej[1][1=\mathbf{j}]{\mathcal{C}_{#1}}   %%%%% typical cone

\newcommandx{\coneindex}[2][2=\mathcal{C}]{#1_#2}

\newcommand\Nset{\mathbb{N}}
\newcommand\Zset{\mathbb{Z}}

\newcommand\Rset{\mathbb{R}}

\newcommandx\borel[1][1=\csms]{\mathcal{B}(#1)}   %% the borel sigma-field of $\Rset^{#1}$

\newcommand{\bszero}{{\boldsymbol0}}

\newcommand\indep{\dag}  %%% notation for independent blocks in beta-mixing blocking method; * was bad; \dag is ugly
\newcommand\bsG{\boldsymbol{G}}

\newcommand\bss{\boldsymbol{s}}

\newcommand\bsx{\boldsymbol{x}}
\newcommand\bsX{\boldsymbol{X}}

\newcommand\bsy{\boldsymbol{y}}
\newcommand\bsY{\boldsymbol{Y}}

\newcommand\bsZ{\boldsymbol{Z}}

\newcommand\bsTheta{\boldsymbol{\Theta}}

\newcommand\bsnu{\boldsymbol{\nu}}

\newcommand\tailmeasure{\bsnu} %%%% the tail measure of a stationary time series
\newcommand\tailmeasurestar{{\bsnu}^*} %%%% the "center" of the tail measure of a stationary time series in terms of $Q$
\newcommand\exc{\mce}  %%%% number of exceedences over 1
\newcommand\clusterlength{\mcl}
\newcommand{\extremalindexfunc}{\mct}
\newcommand\anchor{\mca} %%%%%%% generic notation for anchoring sequences
\newcommandx{\numult}[2][2=]{\bsnu_{\boldsymbol{#1}_{#2}}} %%%% same with an optional argument for index of vector
\newcommand\backest\Upsilon

\newcommandx{\nualphap}[2][1=\alpha,2=p]{\nu_{#1,#2}}  %%%% mesure of regular variation in dimension 1
\newcommandx\numultcondi[2][2=]{\boldsymbol{\nu}_{\boldsymbol{#1}_{#2}}}

\newcommandx\chain[1][1=Y]{\mathbb{#1}}
\newcommandx{\scalingseq}[1][1=n]{c_{#1}}  %Scaling in multivariate RV
\newcommandx{\scalingfunction}[1][1=]{c#1} %Scaling in multivariate RV

 \newcommandx{\absquantileseq}[2][1=]{a_{#1#2}}
\newcommandx{\dhinterseq}[1][1=n]{r_{#1}}  % the intermediate sequence in condition dh
\newcommandx{\dhinterseqsmall}[1][1=n]{\ell_{#1}}  % the intermediate small block sequence (related to condition dh)
\newcommandx{\dinterseq}[1][1=n]{r_{#1}}  % the intermediate sequence in condition d
\newcommandx{\tepseq}[1][1=n]{u_{#1}}  % the sequence $u_n$ to define the TEP

\newcommandx{\interseq}[1][1=n]{k_{#1}}  % intermediate sequence in TEP

\newcommandx{\scalingseqcone}[1][1=n]{c_{#1}} %Scaling sequence for MRV on cone

\newcommandx{\scalingseqhidden}[1][1=n]{\tilde{c}_{#1}}
\newcommandx{\scalingfunctionhidden}[1][1=]{\tilde{c}{#1}}
\newcommandx{\scalingseqcev}[1][1=n]{c^*_{#1}}

\newcommand\conditionS{\ensuremath{\mathcal{S}}} % condition for mean cluster size and for covariance of ordinary TEP
\newcommandx{\conditiondh}[2][1=\dhinterseq,2=\scalingseq]{\ensuremath{\mathcal{A}\mathcal{C}(#1,#2)}} % condition for tail process to tend to 0
\newcommandx{\conditionANSJB}[1][1=\scalingseq]{\mathrm{ANSJB}(\dhinterseq,#1)}
\newcommandx{\conditiondhsum}[1][1=\scalingseq]{\ensuremath{\mathcal{S}(\dhinterseq,#1)}} % condition for convergence of covariance of TEP
\newcommandx{\conditiondhsumW}[1][1=\scalingseq]{\ensuremath{\mathcal{SW}(\dhinterseq,#1)}} % condition for convergence of covariance of TEP

\newcommand\conditionR{\ensuremath{\mathcal{R}}} % $r_n\bar{F}(u_n)\to0$
\newcommand\convfidi{\stackrel{\mbox{\tiny fi.di.}}{\longrightarrow}} %%% weak convergence of fidis of stochastic processes
\newcommand\convdistr{\stackrel{\mbox{\tiny\rm d}}{\longrightarrow}} % weak convergence of random variables
\newcommand\convprob{\stackrel{\tiny \mathbb{P}}{\longrightarrow}}

\newcommandx\prohodistsym[1][1=]{\rho_{#1}}
\newcommandx\prohodist[3][3=]{\rho_{#3}(#1,#2)}

\newcommand\rmd{\mathrm{d}} % Variable in integrals

\newcommand\esp{\mathbb E}
\newcommand\pr{\mathbb P}
\newcommand\var{\mathrm{Var}}
\newcommand\cov{\mathrm{Cov}}

\newcommandx{\autocov}[1][1=]{\gamma_{#1}}

\newcommandx{\cdfnorm}[1][1=\bsX]{H_{#1}}
\newcommandx{\tailcdfnorm}[1][1=\bsX]{\overline{H}_{#1}}

\newcommand\ind[1]{\mathbbm{1}{\left\{#1\right\}}}

\newcommand\mca{\mathcal A}

\newcommand\mce{\mathcal E}
\newcommand\mcf{\mathcal F}
\newcommand\mcg{\mathcal G}

\newcommand\mch{\mathcal H}

\newcommand\mcl{\mathcal L}

\newcommand\mct{\mathcal T}

\newcommandx\test[2][1=X]{{#1}_{#2}}
\newcommandx\orderstat[3][1=X]{{#1}_{(#2:#3)}}
\newcommand\statinterseq{k}

\newcommandx{\sequence}[3][2=\Zset,3=j]{\{#1_{#3},#3\in#2\}}
\newcommandx{\sequenceshort}[2][2=j]{\{#1_#2\}}
\newcommandx\sequ[3][2=j,3=\mathbb{Z}]{\{#1_#2,#2\in#3\}}
\newcommandx\sequnorm[3][3=j,2=\mathbb{Z}]{\{\norm{#1_#3},#3\in#2\}}
\newcommandx\sequnormq[4][2=,4=j,3=\mathbb{Z}]{\{\norm{#1_#4}^{#2},#4\in#3\}}

\newcommandx\uncompactd[2][1=d]{(\overline{\Rset}^{#1})^{#2}\setminus\{\boldsymbol0\}}
\newcommandx{\barrsetproduct}[2][1=d]{(\overline{\Rset}^{#1})^{#2}}
\newcommandx{\rsetproduct}[2][1=d]{(\Rset^{#1})^{#2}}

\newcommand{\metricspace}{\csms}
\newcommandx\csms[1][1=E]{\mathsf{#1}}   %%% a generic complete separable metric space
\newcommandx\borelcsms[1][1=E]{\mathcal{#1}}   %%% the borel sigma-field of the generic complete separable metric space
\newcommandx\mplusx[1][1=]{\mathcal{M}#1}  %%%% the borel measures
\newcommandx\mplusxp[1][1=]{\mathcal{N}{#1}} %%% point measures : with an $N$ as in kallenberg 17
\newcommandx\mplusxpb[1][1=\borelcsms]{\mathcal{N}_{pb}({#1})}  %%% point measures with bounded points
\newcommandx\mplusxpone[1][1=\borelcsms]{\mathcal{N}_{p1}({#1})}  %%% point measures with largest point with modulus 1
\newcommandx\mplusxpeps[1][1=\borelcsms]{\mathcal{N}_{p\epsilon}({#1})}  %%% point measures with one point with modulus greater than $\epsilon$.
\newcommandx\mplusxf[1][1=\borelcsms]{\mathcal{M}_f({#1})}
\newcommandx\mplusxpf[1][1=\borelcsms]{\mathcal{N}_{pf}({#1})} %%%% finite point measures
\newcommandx\mplusxps[1][1=\borelcsms]{\mathcal{N}_{ps}({#1})} %%%% simple point measures
\newcommandx\mplusxpS[1][1=\borelcsms]{\mathcal{N}_{pS}({#1})} %%%%  point measures  with summable points
\newcommandx\mplusxpsc[1][1=\borelcsms]{\mathcal{N}_{psc}({#1})} %%%%  point measures with jumps of constant sign

\newcommand\distance{\mathrm{d}}  %%%% generic distance
\newcommandx\metric[1][1=\metricspace]{\distance_{#1}} %%%% metric for a general metric space $\metricspace$
\newcommandx\metricmcg{\rho} %%%% metric for the index class $\cg$

\newcommandx\bracknum[3][2=\mch]{N_{[\,]}(#1,#2,#3)} % not needed anymore
\newcommandx\bracknumarray[2][2=\mch]{N_{[\,]}(#1,#2,L^2_n)} % bracketing number for arrays - in particular for tep % not needed anymore
\newcommandx\entropynum[3][3=\mch]{N(#1,#3,#2)} % not needed anymore

\newcommandx\process[1][1=X]{\mathbb{#1}}

\newcommandx\hillest[3][1=n,3=]{\widehat{\gamma}_{#1,#2}^{#3}}
\newcommandx\hillmoment[2][1=n]{\widehat{\gamma}_{#1,#2}^{(M)}}

\newcommand\lzero{\ell_0}

\newcommandx\lalpha[1][1=\alpha]{\ell_{#1}}

\newcommand\iid{i.i.d.}

\newcommand\wrt{with respect to}

\newcommand\nonnegative{non-negative}

\begin{document}
\title{Limit theorems for unbounded cluster functionals of regularly varying time series}

\author{Zaoli Chen\thanks{University of Ottawa} \and Rafa{\l} Kulik\thanks{University of Ottawa, email: rkulik@uottawa.ca}}

\date{\today}
\maketitle

\begin{abstract}
A “blocks method” is used to define clusters of
extreme values in stationary time series. 
The cluster starts at the first large value in the block
and ends at the last one. The block cluster measure (the point measure at clusters) encodes different aspects of extremal properties. 
Its limiting behaviour is handled by vague convergence, hence the set of test functions consists of bounded, shift-invariant functionals that vanish around zero. If unbounded or non shift-invariant functionals are considered, we may obtain convergence at a different rate, depending on the type of the functional and the block size (small vs~.~large blocks). There are two prominent examples of such functionals: the locations of large jumps and the cluster length. We obtain a comprehensive characterization of the limiting behaviour of the block cluster measure evaluated at such functionals for stationary, regularly varying time series. 

Once the convergence of the block cluster measure is established, we can proceed with consistency of the empirical cluster measure. Consistency holds in the small and moderate blocks scenario, while fails in the large blocks situation. 
Next, we continue with weak convergence of the empirical cluster processes. The starting point is the seminal paper by Drees and Rootzen (2010). Under the appropriate uniform integrability condition (related to small blocks) the results in the latter paper are still valid. In the moderate and large blocks scenario, the Drees and Rootzen empirical cluster process diverges, but converges weakly when re-normalized properly. 
\end{abstract}
\setcounter{tocdepth}{3}
\tableofcontents

\section{Introduction}\label{sec:intro}
For the last several years, one of the main challenges in extreme value statistics has been modelling and estimation of clusters of extremes values. These clusters stem from temporal dependence in time series that represent real-world phenomena. Cluster indices can be used to quantify different aspects of clustering of large values. The classical extremal index, the quantity that appears as an additional parameter in the limiting distribution of maxima in stationary sequences, is one example of the cluster index. The extremal index achieves values in $(0,1]$ such that smaller values indicate stronger clustering.  

Cluster indices are estimated by \textit{blocks estimators}. Early results focused on asymptotic normality of different types of estimators of the extremal index including, among others,  \cite{smith:weissman:1994}, \cite{hsing:1991}, \cite{hsing:1993}. Asymptotic results for tail array sums (see \cite{rootzen:leadbetter:dehaan:1998}) as well as for cluster size distribution (see \cite{robert:2009}) followed. We note that these papers deal primarily with Peak-over-Threshold (PoT) framework, as opposed to the block maxima (BM) method (see \cite{bucher:zhou:2018} for a review of both PoT and BM approaches). The PoT framework is employed in the paper. 

The seminal paper \cite{drees:rootzen:2010} provides an unified approach for obtaining limiting results for clusters of extremes in stationary
sequences. The mathematical set-up of \cite{drees:rootzen:2010} (adapted to the needs of this paper, in the spirit of \cite{kulik:soulier:2020}) is as follows. We consider 
a stationary, regularly varying $\Rset^d$-valued time series $\bsX=\sequence{\bsX}$. 
Let $\dhinterseq$ be a sequence of nonnegative integers diverging to infinity slower than $n$. The numbers $\dhinterseq$ play a role of the block size, needed to capture clustering of large values. Thus, we have $m_n:=[n/\dhinterseq]=\max\{j\in \mathbb{N}:j\leq n/\dhinterseq\}$ blocks. Let $0<\tepseq\to\infty$ and set $w_n:=\pr(\norm{\bsX_0}>\tepseq)$, where $|\cdot|$ is a norm on $\Rset^d$. Then 
\begin{align}
\label{eq:random-element-entire}
\mathbb{X}_{j}:=\tepseq^{-1}\left(\bsX_{(j-1)\dhinterseq+1},\ldots,\bsX_{j\dhinterseq}\right)\;, \ \ j=1,\ldots,m_n\;,
\end{align}
is the (scaled) $j$th block and 
\begin{align*}
\tailmeasurestar_{n,\dhinterseq}:=\frac{1}{\dhinterseq w_n}\esp[\delta_{\mathbb{X}_1}]
\end{align*}
is the block cluster measure; see Eq. (6.2.1) in \cite{kulik:soulier:2020}.

For a given functional $H$ in a suitable function class $\mch$, the cluster index is obtained by evaluating the block cluster measure at $H$. It is defined as the limit
\begin{align}
\tailmeasurestar(H):=\lim_{n\to\infty}\tailmeasurestar_{n,\dhinterseq}(H)=
\lim_{n\to\infty} \frac{\esp[H(\mathbb{X}_1)]}{\dhinterseq w_n}=\lim_{n\to\infty} \frac{\esp[H(\tepseq^{-1}\bsX_{1,\dhinterseq})]}{\dhinterseq w_n}\;, \ \ H\in\mch\;.  \label{eq:thelimitwhichisnolongercalledbH}
\end{align}
See \cite[Chapter 6]{kulik:soulier:2020}. The rate of convergence $\dhinterseq w_n$
is asymptotically proportional to the probability of an occurrence of a large value in a single block. The "Peak-over-Threshold (POT)" condition (see \ref{eq:rnbarFun0} below) yields that $\dhinterseq w_n\to 0$. That is, with growing $n$, we are less and less likely to observe exceedances over the high threshold $\tepseq$.

Functionals $H$ are supposed to capture extremal structure of the underlying time series. Thus, they are defined on $(\Rset^d)^\Zset$ and are such that their values do not depend on coordinates whose entries are small.
More precisely, for $\bsX=\{\bsX_j,j\in\Zset\}\in (\Rset^d)^\Zset$ and $i\leqslant j\in\Zset$,  we denote $\bsX_{i,j}=(\bsX_i,\ldots, \bsX_j)\in (\Rset^d)^{(j-i+1)}$.
Then, we identify $H(\bsX_{i,j})$ with
$H((\bszero,\bsX_{i,j},\bszero))$, where $\bszero\in (\Rset^d)^\Zset$ is the zero sequence. Such functionals $H$ will be called \textit{cluster functionals}. Moreover, we will assume that $H(\bsx)=0$ whenever $\bsx^*:=\sup_{j\in\Zset}|\bsx_j|<\epsilon$ for some $\epsilon>0$. For the latter property, we will say that \textit{$H$ vanishes around zero}. Alternatively, we will say that \textit{$H$ has support separated from zero}. See \Cref{sec:classes} for the precise statement. Some examples of cluster functionals include:
\begin{itemize}
\item $H(\bsx)=\ind{\bsx^*>1}$ with $\bsx\in (\Rset^d)^{\Zset}$ and $\bsx^*=\sup_{j\in\Zset}\norm{\bsx_j}$. This functional leads to the block estimator of the extremal index.
\item $H(\bsx)=\ind{\sum_{j\in\Zset}\ind{\norm{\bsx_j}>1}=m}$ leads to the block estimator of the cluster size distribution (cf. the relevant estimation problem in \cite{robert:2009}). 
\item $H_\phi(\bsx)=\sum_{j\in\Zset}\phi(\bsx_j)$ with $\phi:\Rset^d\to\Rset$ gives a tail array sum (cf. \cite{rootzen:leadbetter:dehaan:1998}). The choice $\phi(\bsy)=\ind{y_1>s_1,\ldots,y_d>s_d}$, $\bsy=(y_1,\ldots,y_d)$, $\bss=(s_1,\ldots,s_d)$, $s_i>0$, leads to the multivariate tail empirical process indexed by $\bss$; cf. Example 3.1 in \cite{drees:rootzen:2010}. 
\end{itemize}

The appropriate tool to study the limit in \eqref{eq:thelimitwhichisnolongercalledbH} is \textit{vague convergence of cluster measures}.  Indeed, 
Theorem 6.2.5 in \cite{kulik:soulier:2020} (see also \Cref{theo:cluster-RV} below) establishes such convergence, where the class $\mch$ of test functions consists of continuous and \textit{bounded}, \textit{shift-invariant} functionals that vanish around zero.

\paragraph{Unbounded cluster functionals.}
In this paper we want to study general functionals. Besides the natural mathematical generalisation, there are several motivations for this.   
A cluster length, the distance between the location of the last and the first large value in a block, plays an important role in \cite{drees:rootzen:2010}. Furthermore, it plays a very special role in the context of asymptotic expansion of blocks estimators, the problem studied in \cite{chen:kulik:2023b}. Note that the cluster length is shift-invariant, but unbounded. Hence, at the first step of our analysis we need to extend vague convergence of the block cluster measure to convergence of unbounded cluster functionals.  Their growth will be controlled by a power of the cluster length and as such we will refer to them as the \textit{cluster length-type} functionals. It turns out that such an extension is non-trivial. Under the appropriate uniform integrability condition we can recover \eqref{eq:thelimitwhichisnolongercalledbH} - cluster functionals still have the rate $\dhinterseq w_n$. The uniform integrability condition holds as long as the new anticlustering condition is satisfied. 
The latter in turn is related to \textit{small blocks}. On the other hand, 
in case of \textit{large blocks}, \eqref{eq:thelimitwhichisnolongercalledbH} is no longer true and cluster functionals converge at a different rate. 

\paragraph{Beyond shift-invariance.}
Likewise, for some unbounded non shift-invariant functionals, \eqref{eq:thelimitwhichisnolongercalledbH} does not hold for any choice of the blocks size. A particular example here is the location of large values in a stationary time series.  

In summary,  in the first step of our analysis we study the limit 
\begin{align*}
\lim_{n\to\infty}\frac{\esp[H(\tepseq^{-1}\bsX_{1,\dhinterseq})]}{a_n}\;,
\end{align*}
for appropriately chosen functional $H$ and a suitable sequence $a_n\to 0$. 
The main results that pertain to this problem are as follows. 
\begin{itemize}
\item Distribution of the locations of jumps in a single block. They are uniformly distributed and collapse to one random point. The limiting behaviour does not depend on the block size. See \Cref{sec:single-block-jumps}.
\item Extension of vague convergence of clusters (Theorem 6.2.5 in \cite{kulik:soulier:2020}) to unbounded functionals under uniform integrability conditions and phase transition with respect to block size. 
    See \Cref{sec:extentions-cluster-length}. 
\end{itemize}

\paragraph{Consistency of the empirical cluster measure.}
An empirical counterpart of $\tailmeasurestar_{n,\dhinterseq}$ is the \textit{empirical cluster measures}
\begin{align*}
\widetilde{\bsnu}_{n,\dhinterseq}^*=\frac{1}{nw_n}\sum_{j=1}^{m_n}\delta_{\mathbb{X}_j}\;;
\end{align*}
cf. Eq. (10.0.1) in \cite{kulik:soulier:2020}. Then, $\widetilde{\bsnu}_{n,\dhinterseq}^*(H)$ can be considered as a disjoint block pseudo-estimator of the cluster index $\tailmeasurestar(H)$. 
When \eqref{eq:thelimitwhichisnolongercalledbH} holds, we can expect that $\widetilde{\bsnu}_{n,\dhinterseq}^*(H)$ is consistent. Indeed, for the cluster length-type functionals, consistency holds under 
\textit{the small and the moderate blocks conditions} (see \eqref{eq:small-blocks-ep}-\eqref{eq:moderate-blocks-ep}). For \textit{the large blocks} (see \eqref{eq:large-blocks-ep}) the pseudo-estimator diverges. With a different normalization, we can achieve convergence in probability, but the limit differs from $\tailmeasurestar(H)$. In other words, in the large blocks scenario, we do not obtain a consistent estimator of the cluster index. 
See Sections \ref{sec:clt-small-blocks}-\ref{sec:clt-large-blocks}. 

Likewise, when $H$ records location of a large jumps, then $\widetilde{\bsnu}_{n,\dhinterseq}^*(H)$ diverges. On the other hand, when $H$ is the product of the "location of a large jump" and a cluster functional $G$, then after the proper re-normalization we can obtain a consistent pseudo-estimator of $\tailmeasurestar(G)$. These type of estimators can help with reducing a downward bias when approximating cluster indices such as the extremal index. See \Cref{sec:jump-locations}. 

\paragraph{Weak convergence of empirical cluster process(es).}  
The empirical cluster process $\widetilde{\mathbb{G}}_n$ is defined as 
\begin{align}\label{eq:cluster-process}
\widetilde{\mathbb{G}}_n(H)=\frac{1}{\sqrt{nw_n}}\sum_{j=1}^{m_n}\left\{H(\mathbb{X}_{j})-\esp[H(\mathbb{X}_{j})]\right\}\;, 
\end{align}
where $H$ belongs to a suitably chosen function class.  
The functional central limit theorem for the empirical cluster process is established in \cite{drees:rootzen:2010}. In the context of regularly varying time series the result is given in \cite{kulik:soulier:2020}; see Theorem 10.2.1 there. 

For the cluster length-type functionals, the empirical cluster process $\widetilde{\mathbb{G}}_n$ converges in the small blocks scenario, while it diverges when moderate blocks are considered. We need to modify the empirical cluster process. Recall that in both scenarios, the pseudo-estimator is consistent. 
See Sections \ref{sec:clt-small-blocks}-\ref{sec:clt-moderate-blocks}. 

In case of functionals that record jump locations, again a modified empirical cluster process has to be introduced. See \Cref{sec:jump-locations}. 

\paragraph{Statistical consequences.}
The take out message of this paper is that the "classical" theory of weak convergence of clusters (as introduced in \cite{drees:rootzen:2010}) may fail in case of unbounded functionals. The asymptotics may depend on the type of the functional and the block size. Large blocks should be avoided, since they may lead to either different central limit theorem or inconsistency of blocks estimators.  

On the other hand, the blocks estimators based on unbounded functionals may serve as alternative way of estimating cluster indices. In particular, they may help with a downward bias when estimating the classical extremal index. See \Cref{xmpl:extremal-index-estimator}.

\paragraph{Block size terminology.} For convergence of cluster functionals, we distinguish between \textit{small} and \textit{large} blocks. The notion small/large is relative; it depends on the type of functional. For consistency and weak convergence, we have to consider \textit{small}, \textit{moderate} and \textit{large} blocks. Again, these notions are relative. See \Cref{sec:block-size}.

\subsection{Structure of the paper}

\Cref{sec:prel} consists preliminaries. It introduces the tail process, the relevant class of functions (including the cluster length, which appears to be the most important functional in the context of the paper; see \Cref{sec:classes}), cluster measures and cluster indices. 
\Cref{sec:anticlustering-condition} introduces different types of anticlustering conditions. \Cref{sec:cluster-measure-convergence} recalls 
vague convergence of cluster measures. \Cref{theo:cluster-RV} is the most relevant result in this context. 

In \Cref{sec:internal-clusters} we state a variety of results for unbounded cluster functionals.
Locations of jumps (\Cref{sec:single-block-jumps}) and 
functionals dominated by the cluster length (\Cref{sec:extentions-cluster-length}) play a special role. In the second case, we have to distinguish between small and large blocks. The entire contents of this section is new. 

\Cref{sec:weak-convergence} deals with consistency of blocks estimators and weak convergence of the corresponding empirical cluster processes. The results there are new, however, their proofs follow in a relatively straightforward way from \cite[Chapter 10]{kulik:soulier:2020}.

\Cref{sec:technical-details,sec:proofs-clt} contain all the proofs for unbounded functionals and weak convergence of the empirical cluster processes, respectively. 

The goal of \Cref{sec:internal-clusters,sec:weak-convergence} is to state clearly the results, without too much additional "noise". 
A variety of additional remarks can be found in \Cref{sec:technical-details}.

\subsection{What is missing?}
We provide a comprehensive theory in the small blocks scenario. The behaviour in the large blocks situation is limited to some specific functionals and under rather restrictive dependence assumptions. This can be viewed as a limitation of the paper, however, one can argue that large blocks scenario is much less relevant from the statistical perspective.
 
Also, we do not consider here functionals that do not vanish around zero, for example those that lead to the large deviation index (see e.g. \cite{mikosch:wintenberger:2016}). These functionals may be large due to a cumulation of small values. This may be prevented by imposing a "negligibility of small values" condition, but it is still not suitable for the techniques used in the paper.

\section{Preliminaries}\label{sec:prel}

In this section we fix the notation and introduce the relevant classes of functions. In \Cref{sec:tail-process} we recall the notion of the tail and the spectral tail process (cf. \cite{basrak:segers:2009}). In \Cref{sec:cluster-index} we define cluster indices; see \cite[Chapter 5]{kulik:soulier:2020} for a detailed introduction.

In \Cref{sec:anticlustering-condition} we introduce  anticlustering conditions. The first one, the classical \ref{eq:conditiondh} is needed, to establish convergence of the cluster measure; see
\Cref{sec:cluster-measure-convergence}.  The results of the latter section are extracted from \cite[Chapter 6]{kulik:soulier:2020}.
See also \cite{planinic:soulier:2018} and \cite{basrak:planinic:soulier:2018}. Another condition, called here \ref{eq:conditionSstronger:gamma}, is crucial for many results in the paper.

In \Cref{sec:dependence-assumptions} we introduce our dependence assumptions.

\subsection{Notation}
Let $\norm{\cdot}$ be an arbitrary norm on $\Rset^d$ and $\{\tepseq\}$, $\{\dhinterseq\}$ be such that
\begin{align}\label{eq:rnbarFun0}
\lim_{n\to\infty}\tepseq=\lim_{n\to\infty}\dhinterseq =\lim_{n\to\infty}nw_n = \infty\;, \ \ \lim_{n\to\infty}\dhinterseq/n=\lim_{n\to\infty}\dhinterseq w_n = 0\;,
\tag{$\conditionR(\dhinterseq,\tepseq)$}
\end{align}
where $w_n=\pr(\norm{\bsX_0}>\tepseq)$. 
For a sequence $\bsx=\{\bsx_j,j\in\Zset\}\in (\Rset^d)^\Zset$ and $i\leqslant j\in\Zset\cup\{-\infty,\infty\}$ denote $\bsx_{i,j}=(\bsx_i,\ldots,\bsx_j)\in (\Rset^d)^{j-i+1}$, $\bsx_{i,j}^\ast=\max_{i\leqslant l\leqslant j}|\bsx_l|$ and $\bsx^\ast=\sup_{j\in\Zset}|\bsx_j|$. By $\bszero$ we denote the zero sequence; its dimension can be different at each of its occurrences.

By $\lzero(\Rset^d)$ we denote the set of $\Rset^d$-valued sequences which tend to zero at infinity.

\subsection{Tail process}\label{sec:tail-process}
Let $\bsX=\sequence{\bsX}$ be a stationary, regularly varying time series with values in $\Rset^d$ and tail index $\alpha$. 
 Then, there exists a sequence $\bsY=\sequence{\bsY}$ such that
\begin{align*}
  %\label{eq:def-numultprob}
  \pr(x^{-1}(\bsX_i,\dots,\bsX_j) \in \cdot \mid |\bsX_0|>x)   \mbox{ converges weakly to }  \pr((\bsY_i,\dots,\bsY_j) \in \cdot)
\end{align*}
as $x\to\infty$ for all $i\leqslant j\in\Zset$.
We call $\bsY$ the tail process.
See \cite{basrak:segers:2009}. Equivalently, viewing $\bsX$ and $\bsY$ as random elements with values in $(\Rset^d)^\Zset$, we have for every bounded or
\nonnegative\ functional~$H$ on $(\Rset^d)^{\Zset}$, continuous \wrt\ the law of the tail process $\bsY$,
\begin{align}\label{eq:tailprocesstozero}
   \lim_{x\to\infty} \frac{\esp[H(x^{-1}\bsX)\ind{\norm{\bsX_0}>x}]} {\pr(\norm{\bsX_0}>x)}
  & =  \esp[H(\bsY)] \; .
\end{align}
The random variable $|\bsY_0|$ has the Pareto distribution with index $\alpha$ and hence $\norm{\bsY_0}>1$.
The spectral tail process $\bsTheta=\{\bsTheta_j,j\in\Zset\}$ is defined as $\bsTheta_j=\bsY_j/|\bsY_0|$.

\subsection{Classes of unbounded functions}\label{sec:classes}
Let $\epsilon>0$. 
Define 
\begin{subequations}
\begin{align}
T^{(1)}(\bsx,\epsilon)=T_{\rm min}(\bsx,\epsilon)  &= \inf\{ j \in \Zset : | \bsx_j | > \epsilon \}\;, \label{eq:exc-times-x-1}\\
     T_{\rm max}(\bsx,\epsilon)&=\sup\{ j \in \Zset : | \bsx_j | > \epsilon \}\;,  \label{eq:exc-times-x-2}\\
     T^{(i+1)}(\bsx,\epsilon)&=\inf\{j>T^{(i)}(\bsx,\epsilon): |\bsx_j|>\epsilon\}\;, \ \ i\geqslant 1\;.\label{eq:exc-times-x-3}
\end{align}
\end{subequations}
We will use the convention $\inf\{\emptyset\}=+\infty$, $\sup\{\emptyset\}=-\infty$. If 
$\bsx\in \ell_0(\mathbb{R}^d)$, then $\sup\{ j\geqslant 0: | \bsx_j | > \epsilon   \}<\infty$ and
$\inf\{ j\geqslant 0: | \bsx_j | > \epsilon   \}>-\infty$. Hence, 
when restricted to $\ell_0$, the functionals $T^{(i)}$, $i\geqslant 1$, and $T_{\rm max}$ attain finite values. 

Formally, the above functionals depend on $\epsilon$. However, without loss of generality we will assume throughout the paper that $\epsilon=1$ and we will use the notation $T^{(i)}(\bsx)$ instead. 

Let again $\bsx\in \ell_0(\mathbb{R}^d)$. We can define the
\textbf{cluster length functional} as
\begin{align}\label{eq:cluster-length-def}
\clusterlength(\bsx) =T_{\rm max}(\bsx)-T_{\rm min}(\bsx)+1\;
\end{align}
 with the convention $\clusterlength(\bsx)=0$ whenever $\bsx^\ast=\sup_{j\in\Zset}\norm{\bsx_j}\leqslant 1$. We note also that if $\norm{\bsx_0}>1$, while $\sup_{j\in\Zset, j\not=0}\norm{\bsx_j}\leqslant 1$, then $\clusterlength(\bsx)=1$. 

We define further  
the functional $\exc$ by $\exc(\bsx)=\sum_{j\in\Zset}\ind{\norm{\bsx_j}>1}$. It returns the number of exceedances over 1. 

We shall consider functionals under the following assumptions.
\begin{hypothesis}
\label{Assumption:class-mathcalH}
We denote by $\mch$ a collection of functionals $\ell_0(\Rset^d) \to \Rset_+$ such that  each $H \in \mathcal{H}$ satisfies:
\begin{itemize}
    \item[$(\rom1)$] $H$ is continuous with respect to the law of the process $\bsY$;
    \item[$(\rom2)$] If $\exc(\bsx)=0$, then $H(\bsx)=0$;
    \item[$(\rom3)$] If $\exc(\bsx)>0$, then $H(\bsx) =H(\bsx_{T_{\rm min}(\bsx),T_{\rm max}(\bsx)})$, where
    $T_{\rm min}(\bsx)$ and $T_{\rm max}(\bsx)$, are the first and the last exceedance times defined in \eqref{eq:exc-times-x-1}-\eqref{eq:exc-times-x-2}.
\end{itemize}
\end{hypothesis}
We note that the above assumption allows for unbounded functionals. We will need to control their growth.
\begin{hypothesis}
\label{Assumption:class-mathcalHcH}
We denote by $\mch(\gamma)\subseteq \mch$ a collection of functionals $\ell_0(\Rset^d) \to \Rset$ such that  each $H \in\mch(\gamma)$ satisfies:
\begin{itemize}
\item[$(\rom4)$] There exists a constant $C_H>0$ such that
    $ H(\bsx)  \leqslant C_H  \big[ \clusterlength(\bsx) \big]^{\gamma}$ for all $\bsx \in \ell_0( \Rset^d )$.
\end{itemize}
\end{hypothesis}
\begin{example}{\rm 
Obviously, we can take $H$ as the cluster length functional itself: $H=\clusterlength$. Then $\gamma=1$.

The extremal index functional defined as $\extremalindexfunc(\bsx)=\ind{\bsx^\ast>1}$. It fulfills \Cref{Assumption:class-mathcalHcH} with $C_{\extremalindexfunc}=1$ and $\gamma=0$. This functional is bounded. In particular, $\extremalindexfunc(\bsY)=1$. 

In case of extremal independence, $\norm{\bsY_j}=0$ for $|j|\geqslant 1$ and hence $\extremalindexfunc(\bsY)=\clusterlength(\bsY)=1$.
}\end{example}
\begin{remark}{\rm
Functional can depend on small values, but only those that occur between large values ("within a cluster"). For example, $H_{\rm sum}(\bsx)=\ind{\sum_{j=T_{\rm min}(\bsx)}^{T_{\rm max}(\bsx)}|\bsx_j|>\eta}$ (with $\eta>0$) fulfills \Cref{Assumption:class-mathcalH}, but $H(\bsx)=\ind{\sum_{j\in\Zset}|\bsx_j|>1}$ does not. The latter functional appears in the context of large deviations, see \cite{mikosch:wintenberger:2016}.

We also note that the functionals $T^{(i)}$ defined in \eqref{eq:exc-times-x-3} do not fulfill \Cref{Assumption:class-mathcalH}(iii). Furthermore, they cannot be bounded by the cluster length, and hence they do not fulfill \Cref{Assumption:class-mathcalHcH}.  
}\end{remark}
The class $\mch(\gamma)$ is parametrized by $\gamma$ that will play an important role. The case $\gamma=0$ corresponds to bounded functionals. If a function $H$ is bounded, we will denote $\|H\|=\sup_{\bsx\in(\Rset^d)^\Zset}|H(\bsx)|$. The property $(\rom2)$ will be referred to as "$H$ vanishes around $\bszero$" or "$H$ has support separated from $\bszero$". In particular, $\clusterlength$, $\extremalindexfunc$ and the aforementioned $H_{\rm sum}$ vanish around $\bszero$. However, the large deviation functional does not vanish around $\bszero$.

\subsection{Cluster measure and cluster indices}\label{sec:cluster-index}
Consider
the infargmax functional $\anchor_0$ defined on $(\Rset^d)^\Zset$ by
    $\anchor_0(\bsx)=\inf\{j:\bsx_{-\infty,{j}}^\ast=\bsx^\ast\}$, with the convention that
    $\inf\{\emptyset\}=+\infty$.
If  $\pr(\anchor_0(\bsY)\notin\Zset)=0$ then we can define
\begin{align*}
  %\label{eq:canditheta-anchor}
  \canditheta = \pr(\anchor_0(\bsY)=0) \; .
\end{align*}
In fact, $\anchor_0$ can be replaced with any anchoring map (see \cite[Theorem 5.4.2]{kulik:soulier:2020}). Hence, 
\begin{align}\label{eq:canditheta-anchor-conclusion}
\canditheta=\pr(\anchor_0(\bsY)=0)=\pr\left(\bsY_{-\infty,-1}^\ast\leqslant 1\right)=\pr\left(\bsY_{1,\infty}^\ast\leqslant 1\right)\;.
\end{align}
Therefore, $\canditheta$ can be recognized as the (candidate) extremal index. It becomes the usual extremal index under additional mixing and anticlustering conditions.

\begin{definition}[Cluster measure]
  \label{def:clustermeasure}
  Let $\bsY$ and $\bsTheta$ be the tail process and the spectral tail process, respectively, such that $\pr(\lim_{|j|\to\infty} |\bsY_j|=0)=1$. The {cluster measure} is the measure $\tailmeasurestar$ on
  $\lzero(\Rset^d)$ defined by
  \begin{align*}
    %\label{eq:def-tailmeasurestar-premiere}
\tailmeasurestar    = \canditheta \int_0^\infty \esp[\delta_{r\bsTheta}\ind{\anchor_0(\bsTheta)=0}] \alpha r^{-\alpha-1} \rmd r   \; .
\end{align*}
\end{definition}
The measure $\tailmeasurestar$ is
boundedly finite on $(\Rset^d)^\Zset\setminus\{\bszero\}$, puts no mass at $\bszero$ and is
$\alpha$-homogeneous. For every bounded or non-negative functional $H$ such that $H(\bsx)=0$ if
$\bsx^\ast\leqslant 1$ we have
\begin{align}
    \label{eq:cluster-measure}
  \tailmeasurestar(H) &=  \esp[H(\bsY) \ind{\anchor_0(\bsY)=0}]
  =\esp[H(\bsY) \ind{\bsY_{-\infty,-1}^\ast\leqslant 1}]
   \; .
\end{align}
\begin{definition}[Cluster index]
\label{def:cluster-index}
We will call $\tailmeasurestar(H)$ the cluster index associated to the
functional $H$.
\end{definition}

\subsection{Change of measure}\label{sec:change-of-measure}
It is important to notice the presence of $\ind{\bsY_{-\infty,-1}^\ast\leqslant 1}$ in the definition of $\tailmeasurestar$.  
It will be convenient to express some formulas in the language of an auxiliary process $\bsZ$. Assume that $\pr(\lim_{|j|\to \infty} |\bsY_j| = 0) = 1$ and define the process $\bsZ$
as $\bsY$ conditioned on the first exceedance over 1 to happen at time zero, that is
$\bsY_{-\infty,-1}^*\leqslant 1$. Then, combining \eqref{eq:canditheta-anchor-conclusion} with \eqref{eq:cluster-measure}, we obtain 
\begin{align*}%\label{eq:Palm}
\tailmeasurestar(H)=\esp[H(\bsY) \ind{\bsY_{-\infty,-1}^\ast\leqslant 1}]=\canditheta\esp[H(\bsZ)]\;.
\end{align*}
Informally speaking, the distribution of $\bsZ$ is a Palm version of the distribution of $\bsY$.  
See \cite{planinic:2023} and \cite{last:2023} for a relation between the tail process and Palm theory. See also Exercise 5.29 in \cite{kulik:soulier:2020}.

\subsection{Anticlustering conditions}\label{sec:anticlustering-condition}
For each fixed $r\in\Nset$, the distribution of
$\tepseq^{-1}\bsX_{-r,r}$ conditionally on $\norm{\bsX_0}>\tepseq$ converges weakly to the
distribution of $\bsY_{-r,r}$. In order to let $r$ tend to infinity, we must embed all these finite
vectors into one space of sequences. By adding zeroes on each side of the vectors
$\tepseq^{-1}\bsX_{-r,r}$ and $\bsY_{-r,r}$ we identify them with elements of the space
$\lzero(\Rset^d)$. Then $\bsY_{-r,r}$
converges (as $r\to\infty$) to $\bsY$ in $\lzero(\Rset^d)$ if (and only if) $\bsY\in\lzero(\Rset^d)$ almost surely.

However, this is not enough for statistical purposes and we consider the following definition that controls persistence of large values on one block. 
\begin{definition}
  [\cite{davis:hsing:1995}, Condition~2.8]\label{def:DH}
  Condition~\ref{eq:conditiondh}
  holds if for all $s,t>0$,
  \begin{align}
    \label{eq:conditiondh}
    \lim_{\ell\to\infty} \limsup_{n\to\infty}\pr\left(\max_{\ell\leqslant |j|\leqslant  \dhinterseq}|\bsX_j|
    > \tepseq s\mid |\bsX_0|> \tepseq t \right)=0 \; .
    \tag{$\conditiondh[\dhinterseq][\tepseq]$}
  \end{align}
\end{definition}
Condition \ref{eq:conditiondh} is referred to as the (basic) anticlustering condition. 
\ref{eq:conditiondh} implies that $\bsY\in \lzero(\Rset^d)$
and
    $\vartheta = \pr\left ( \bsY^\ast_{-\infty,-1} \leqslant 1  \right)> 0$.
Also, \ref{eq:conditiondh} holds for sequence of \iid\ random
variables whenever $\lim_{n\to\infty} \dhinterseq w_n=0$.

\begin{definition}\label{def:conditionSstronger:gamma}
\hypertarget{SummabilityAC}{Condition
  \ref{eq:conditionSstronger:gamma}} holds if
for all $s,t>0$
\begin{align}
    \label{eq:conditionSstronger:gamma}
    \lim_{\ell\to\infty} \limsup_{n\to\infty} \frac{1}{ w_n}
    \sum_{i=\ell}^{\dhinterseq}i^\gamma \pr(\norm{\bsX_0}>\tepseq s,\norm{\bsX_i}>\tepseq t) = 0 \;
    . \tag{$\conditiondhsumstrongergammapaper{\dhinterseq}{\tepseq}{\gamma}$}
  \end{align}
\end{definition}
This condition is new.
It is an extension of the existing condition \ref{eq:conditionSstronger:gamma} with $\gamma=0$. The latter was proven to hold for a large class of time series models, such as geometrically ergodic Markov chains,
short-memory linear or max-stable processes. See \cite{kulik:soulier:2020}.

\begin{remark}\label{rem:rn2gamma}
{\rm
\begin{itemize}
\item The condition implies that $\sum_{i\in\Zset}i^\gamma \pr(\norm{\bsY_i}>1)<\infty$. See Lemma 9.2.3 in \cite{kulik:soulier:2020} in case $\gamma=0$. Its proof can be easily modified to arbitrary $\gamma>0$.
\item
It is obvious that in case of \iid\ or $\ell$-dependent sequences \ref{eq:conditionSstronger:gamma} holds if and only if $\lim_{n\to\infty}\dhinterseq^{\gamma+1} w_n=0$. However, we could not establish that \ref{eq:conditionSstronger:gamma} and $\lim_{n\to\infty}\dhinterseq^{\gamma+1} w_n=0$ are equivalent in general. 

\item
The condition~\hyperlink{SummabilityAC}{$\mathcal{S}^{(\gamma_1)}(\dhinterseq, \tepseq)$} implies \hyperlink{SummabilityAC}{$\mathcal{S}^{(\gamma_2)}(\dhinterseq, \tepseq)$} whenever $\gamma_1>\gamma_2\geqslant 0$.
    The condition~\hyperlink{SummabilityAC}{$\mathcal{S}^{(0)}(\dhinterseq, \tepseq)$} implies
    the condition~\ref{eq:conditiondh}.
\end{itemize}
}
\end{remark}

\subsection{Conditional convergence of clusters}
From the definition of the tail process we obtain immediately
\begin{align*}
    \lim_{n\to\infty} \esp[H(\tepseq^{-1} \bsX_{i,j})\mid \norm{\bsX_0}>\tepseq] = \esp[H(\bsY_{i,j})] \; ,
    %\label{eq:conv-H-cluster}
  \end{align*}
  for $i\leqslant j$ and suitable functionals $H$. Thanks to the anticlustering condition we can replace $i,j$ with a sequence diverging to $\infty$.

  \begin{proposition}[\cite{basrak:segers:2009}, Proposition~4.2; \cite{kulik:soulier:2020}, Theorem 6.1.4]
  \label{lem:tailprocesstozero}
  Let $\bsX$ be a stationary time series. Assume that~\ref{eq:conditiondh} holds.
  Let $H$ be a bounded or non-negative functional defined on $\ell_0(\Rset^d)$ that is almost
surely continuous with respect to the distribution of the tail process $\bsY$.
  Then for any constant $C>0$,
  \begin{align*}
    \lim_{n\to\infty} \esp[H(\tepseq^{-1} \bsX_{-C\dhinterseq,C\dhinterseq})\mid \norm{\bsX_0}>\tepseq] = \esp[H(\bsY)] \; .
    %\label{eq:conv-H-cluster}
  \end{align*}
\end{proposition}

\subsection{Vague convergence of cluster measure}\label{sec:cluster-measure-convergence}
We now investigate the unconditional convergence of
$\tepseq^{-1}\bsX_{1,\dhinterseq}$. Contrary to
\Cref{lem:tailprocesstozero}, where an extreme value was imposed at time 0, a large value in the
cluster can happen at any time.

Define the measures $\tailmeasurestar_{n,\dhinterseq}$, $n\geqslant1$, on $\lzero(\Rset^d)$ as follows:
\begin{align*}
 % \label{eq:def-tailmeasuresta\dhinterseq}
  \tailmeasurestar_{n,\dhinterseq}
  & =  \frac{1} {\dhinterseq  w_n} \esp \left[\delta_{\tepseq^{-1}\bsX_{1,\dhinterseq}} \right] \; .
\end{align*}
We are interested in convergence of $\tailmeasurestar_{n,\dhinterseq}$ to
$\tailmeasurestar$.
We quote Theorem 6.2.5 in \cite{kulik:soulier:2020}.
\begin{proposition}
\label{theo:cluster-RV}
  Let  condition~\ref{eq:conditiondh} hold. Then, for all bounded continuous shift-invariant functionals
$H$ with support separated from $\bszero$ we have
  \begin{align}\label{eq:vague-convergence}
    \lim_{n\to\infty} \tailmeasurestar_{n,\dhinterseq}(H)
    = \lim_{n\to\infty}\frac{\esp[ H(\tepseq^{-1}\bsX_{1,\dhinterseq})]} {\dhinterseq w_n}
    = \tailmeasure^\ast(H) \; .
  \end{align}
\end{proposition}
The immediate consequence is the following limit:
  \begin{align}
  \label{eq:extremal-index}
    \lim_{n\to\infty} \frac{\pr(\bsX_{1,\dhinterseq}^\ast>\tepseq)}{\dhinterseq w_n}
     =\canditheta  \; .
  \end{align}
If $H$ does not fulfill the conditions of \Cref{theo:cluster-RV}, we may have several scenarios, as studied in \Cref{sec:internal-clusters}. 
\subsection{Dependence assumptions}\label{sec:dependence-assumptions}
We shall divide the sample $\{\bsX_1,\ldots, \bsX_n\}$ into $m_n\in\mathbb{N}$ disjoint blocks of same size $\dhinterseq\in\mathbb{N}$.
Without loss of generality we will assume that $n=m_n\dhinterseq$. Dependence within each block ("local dependence") is controlled by the appropriate anticlustering condition. Dependence between the blocks ("global dependence") is controlled by $\beta$-mixing. We will assume that there exists a sequence $\ell_n$ such that the following condition holds: 
  \begin{align}\label{eq:conditionbeta-for-Aprime}
    \lim_{n\to\infty} \frac1{\ell_n} = \lim_{n\to\infty} \frac{\ell_n}{\dhinterseq} =
    \lim_{n\to\infty} \frac{\dhinterseq}n  = \lim_{n\to\infty} \frac{n}{\dhinterseq} \beta_{\ell_n} = 0 \; .\tag{$\beta(\dhinterseq,\ell_n)$}
  \end{align}
We recall the classical covariance inequality for bounded, beta-mixing random variables. 
Let $\beta({\mcf_1},{\mcf_2})$ be the $\beta$-mixing coefficient between two sigma fields.
Then
\begin{align}\label{eq:mixing-inequality}
|\cov(H(Z_1),H(Z_2))|\leq \constant\ \|H\|_{\infty} \|\widetilde{H}\|_{\infty}\beta(\sigma(Z_1),\sigma(Z_2))\;.
\end{align}
In \eqref{eq:mixing-inequality} the constant $\constant$ does not depend on $H,\widetilde{H}$.

\subsection{Notation for blocks}\label{sec:notation-for-blocks}
We consider disjoint blocks of size $\dhinterseq$:
\begin{align*}
I_{j}=\{(j-1)\dhinterseq+1,\ldots,j\dhinterseq\}\;, \ \ j=1,\ldots,m_n\;.
\end{align*}
The subscript $j$ will always indicate the numbering of blocks. 
Denote by
\begin{align*}
N_{j}:=\exc(\tepseq^{-1}\bsX_{(j-1)\dhinterseq+1,j\dhinterseq})=\sum_{i\in I_{j}}\ind{\norm{\bsX_i}>\tepseq}\;, \ \ j=1,\ldots,m_n\;,
\end{align*}
the number of exceedances in the block $I_{j}$.
Set
\begin{align}\label{eq:set-Aj}
A_{j} = \{ \exists i : (j-1)\dhinterseq + 1 \leqslant i \leqslant  j\dhinterseq , \vert \bsX_i \vert > \tepseq  \}=\{\bsX_{(j-1)\dhinterseq + 1,j\dhinterseq}^\ast>\tepseq\}\;.
\end{align}
Recall the notation \eqref{eq:exc-times-x-1}-\eqref{eq:exc-times-x-3}. 
If $N_{j}\not=0$, then denote by 
\begin{align*}
t_j{(i)}=T^{(i)}(\tepseq^{-1}\bsX_{(j-1)\dhinterseq+1,j\dhinterseq})
\end{align*} 
the $i$th exceedance time in the block $j$. Hence, 
\begin{align}\label{eq:exceedance}
t_j{(1)}<\cdots< t_{j}{(N_{j})}\;, \ \ j=1,\ldots,m_n\;,
\end{align}
are 
the exceedance times in the $j$th block. We will use the convention $t_{j}(0)\equiv (j-1)\dhinterseq$ and $t_{j}(N_{j}+1)\equiv j\dhinterseq$. Note that it is possible that
$t_{j}{(N_{j})}=t_{j}(N_{j}+1)=j\dhinterseq$. This happens when the last jump in the block $j$ occurs at the right-end point, $j\dhinterseq$. On the other hand, since the $j$th block starts at $(j-1)\dhinterseq+1$, $t_{j}(1)$ is strictly larger than $t_{j}(0)$.
We will set $t_1(i)=t(i)$ for $i=1,\ldots,N_1$. 

Furthermore, using the notation introduced in \eqref{eq:cluster-length-def},
\begin{align}\label{eq:def-clusterlength}
t_{j}{(N_{j})}-t_{j}{(1)}+1=\clusterlength(\tepseq^{-1}\bsX_{(j-1)\dinterseq+1,j\dhinterseq})=
:\clusterlength_{j}
\end{align}
is the cluster length in the $j$th block.
We use the convention $\sum_{i=\ell}^j=0$ whenever $j<\ell$.
Hence, $\clusterlength_{j}=1$ whenever $N_{j}=1$. If there are no exceedances over the threshold $\tepseq$, then we set $\clusterlength_{j}\equiv 0$ (there will be no issue with this, since $\clusterlength_{j}$ will appear with the appropriate indicator.) Also, by the definition, $\clusterlength_{j}\leqslant \dhinterseq$ uniformly in $j$.

%We need to keep track of the exceedances in each block.
%For this, for each $j=1,\ldots,m_n$ and each $1\leqslant k_1 \leqslant k_2 \leqslant N_{j}$, we define
%\begin{align*}
%%\label{eq:random-elements-Z}
%\mathbb{X}_{j}{(k_1:k_2)} = \tepseq^{-1} \big(\bsX_{ t_{j}{(k_1)}}, \ldots, \bsX_{t_{j}{(k_2)} }  \big)\;.
%\end{align*}
%Thus, in particular, $\mathbb{X}_{j}{(1:N_{j})}$ is a vector of a random size $\clusterlength_{j}$ that consists of all (scaled) exceedances in the block $j$ and all small values between the first and the last exceedance. We use the convention $\mathbb{X}_{j}{(k_1:k_2)}=0$ whenever $k_1>k_2$. Hence, in particular,
%$\mathbb{X}_{j}{(1:N_{j})}=0$ if $N_{j}=0$.

We will also use the following notation for the scaled block $j$:
\begin{align*}
%\label{eq:random-element-entire}
\mathbb{X}_{j}=\tepseq^{-1}\bsX_{(j-1)\dhinterseq+1,j\dhinterseq}\;, \ \ j=1,\ldots,m_n\;.
\end{align*}
\subsubsection{Block size terminology}\label{sec:block-size}
In what follows, when convergence of cluster functionals is considered, the \textit{small} and \textit{large} blocks refers to $\dhinterseq^{\gamma+1}w_n\to 0$ and $\dhinterseq^{\gamma+1}w_n\to \infty$, respectively. In the context of weak convergence of empirical processes, we will distinguish between three different scenarios, referred to as the \textit{small (respectively, moderate or large) blocks condition for empirical processes}:
\begin{subequations}
\begin{align}
\label{eq:small-blocks-ep}
&\dhinterseq^{2\gamma+1}w_n\to 0  \ \ & {\rm (small\  blocks)}\;,\\
\label{eq:moderate-blocks-ep}
&\dhinterseq^{\gamma+1}w_n\to 0\ {\rm and}\ \dhinterseq^{2\gamma+1}w_n\to \infty \ \ & {\rm (moderate\  blocks)}\;,\\
\label{eq:large-blocks-ep}
&\dhinterseq^{\gamma+1}w_n\to \infty  \ \ & {\rm (large\  blocks)}\;.
\end{align}
\end{subequations}
\section{Convergence of unbounded cluster functionals}\label{sec:internal-clusters}
In what follows, $\bsX$ is $\Rset^d$-valued regularly varying time series. Also, \ref{eq:rnbarFun0} is assumed everywhere. 

\Cref{theo:cluster-RV} is valid for bounded, shift-invariant functionals that vanish around zero. In other words, one deals with convergence of the appropriately scaled $\esp[H(\tepseq^{-1}\bsX_{1,\dhinterseq})\ind{A_1}]$, where $A_1=\{\bsX_{1,\dhinterseq}^\ast> \tepseq\}$. The goal is to extend it to unbounded functionals. 

Set $t(i)=t_{1}(i)=T^{(i)}(\tepseq^{-1}\bsX_{1,\dhinterseq})$, $i=1,\ldots,N_{1}$. Hence, $t(1)=T_{\min}(\tepseq^{-1}\bsX_{1,\dhinterseq})$ and $t{(N_1)}=T_{\max}(\tepseq^{-1}\bsX_{1,\dhinterseq})$. 
\subsection{Summary of the results}\label{sec:summary-internal}
\begin{itemize}
\item 
\Cref{sec:single-block-jumps}: 
The functionals $T_{\min}$ and $T_{\max}$ are unbounded and not shift-invariant. Then, $t(1)$ and $t(N_1)$ are not tight under the conditional law (given $A_1$) - they need to be scaled by the block size $\dhinterseq$ and the limiting conditional distribution is uniform. See \Cref{lem:first-and-last-jump}. Next, \Cref{theo:cluster-RV} is no longer valid - the $\gamma$ moment of $t(1)$ or $t(N_1)$ requires the scaling $\dhinterseq^{\gamma+1}w_n$ instead of $\dhinterseq w_n$. See \Cref{lem:fist-jump-function-of}.
\item \Cref{sec:extentions-cluster-length} deals with functionals bounded by the cluster length. Even though the cluster length is a simple function of $t({1})$ and $t(N_1)$, its behaviour changes significantly. First, the cluster length is tight under the conditional law; see \Cref{lem:clusterlength-tight}. Next, 
in the small blocks scenario the uniform integrability holds and \Cref{theo:cluster-RV}  is still valid. It fails in the large blocks scenario. See \Cref{lem:clusterlength-moments} and \Cref{lem:clusterlength-moments-large}. For the $\gamma$ moment of the cluster length, the dichotomous between the small and the large blocks scenario is related to whether $\dhinterseq^{\gamma+1}w_n\to 0$ or $\dhinterseq^{\gamma+1}w_n\to \infty$. 
\end{itemize}
All proofs are given in \Cref{sec:proof-for-internal}.
\subsection{Location of jumps}\label{sec:single-block-jumps} 
First, we analyse the limiting conditional distribution of the first and the last jump times. Note that these functionals are neither bounded nor shift-invariant. When scaled by the block size, both become uniformly distributed. The uniform limiting conditional distribution yields immediately the convergence of conditional moments. The main results are: \Cref{lem:first-and-last-jump}
\Cref{lem:fist-jump-function-of} and \Cref{cor:first-and-last-jump-exp}. No specific dependence structure is needed; only the basic anticlustering condition \ref{eq:conditiondh} is used.
Next, we deal with joint moments of the first and the last jump. The main result is \Cref{lem:moments-first-and-last-jump}. There, we need a slightly stronger anticlustering condition \hyperlink{SummabilityAC}{$\mathcal{S}^{(0)}(\dhinterseq, \tepseq)$}.
The results there are quite intuitive: the first and the last jump time have to be scaled by the block size $\dhinterseq$, while the presence of the indicator of a large jump forces the additional $\dhinterseq w_n$. See \Cref{cor:first-and-last-jump-exp}. It has to be emphasized that the scaling is the same for any choice of the block size. 

In summary, the results of this section show that \Cref{theo:cluster-RV}  is no longer valid for the unbounded, non shift-invariant functionals $T_{\min}$ and $T_{\max}$.

\subsubsection{Distribution of the first and the last jump time}
The first result shows that the scaled jump times $t(1)$ and $t(N_1)$ are approximately uniformly distributed. Only the basic anticlustering condition \ref{eq:conditiondh} is needed. 
\begin{lemma}\label{lem:first-and-last-jump}
Assume that \ref{eq:conditiondh} holds. Let $q\in (0,1)$. 
\begin{enumerate}
\item[{\rm (i)}] We have
\begin{align*}
\lim_{n\to\infty}\frac{1}{\dhinterseq w_n}\pr(\{t(1)>q\dhinterseq\}\cap A_1)=
\lim_{n\to\infty}\frac{1}{\dhinterseq w_n}\pr(\{t(N_1)>q\dhinterseq\}\cap A_1)=(1-q)\canditheta\;.
\end{align*}
\item[{\rm (ii)}] We have
\begin{align*} %\label{eq:joint-convergence-interior-cluster;1}
   \frac{  \mathbb{E}\left[ \ind{ q\dhinterseq <  t(1) <  \dhinterseq }
    G \left(  \tepseq^{-1}\bsX_{1,\dhinterseq}  \right)
    \right] }{\dhinterseq w_n}
     \to
    (1-q)
    \tailmeasurestar(G)\;.
\end{align*}
for any bounded measurable functional $G: (\mathbb{R}^d) ^ \mathbb{Z} \to \mathbb{R}$ that is continuous with respect to the law of $\bsY$, with support separated from $\bszero$.
\end{enumerate}
\end{lemma}
Since $\pr(A_1)\sim \canditheta \dhinterseq w_n$ (see \eqref{eq:extremal-index}), we can rephrase the above result as: 
\begin{align*}\lim_{n\to\infty}\pr(\{t(1)>q\dhinterseq\}\mid A_1)=1-q\;.\end{align*}
\subsubsection{Moments of the first and the last jump time}
The next statements are almost obvious in view of the limiting conditional distribution of $t(1)/\dhinterseq$ being standard uniform.
\begin{lemma}\label{lem:fist-jump-function-of}
Assume that \ref{eq:conditiondh} holds.
Let $f:\Rset_+\to \Rset_+$ be $\gamma$-homogenous, $\gamma>0$. Then
\begin{align*}%\label{eq:first-jump-function-of}
\lim_{n\to\infty}\frac{\esp[f(t(1))G(\tepseq^{-1}\bsX_{1,\dhinterseq})]}{\dhinterseq^{\gamma+1} w_n}=
\lim_{n\to\infty}\frac{\esp[f(t(N_1))G(\tepseq^{-1}\bsX_{1,\dhinterseq})]}{\dhinterseq^{\gamma+1} w_n}=\frac{f(1)}{\gamma+1}\tailmeasurestar(G)
\;
\end{align*}
for any bounded measurable functional $G: (\mathbb{R}^d) ^ \mathbb{Z} \to \mathbb{R}$ that is continuous with respect to the law of $\bsY$, with support separated from $\bszero$. 
\end{lemma}
\begin{corollary}\label{cor:first-and-last-jump-exp}
Assume that \ref{eq:conditiondh} holds. Let $\gamma>0$. Then
\begin{align*}
\lim_{n\to\infty}\frac{\esp[t(1)^\gamma \ind{A_1}]}{\dhinterseq^{\gamma+1} w_n}=\lim_{n\to\infty}\frac{\esp[(\dhinterseq-t(1))^\gamma \ind{A_1}]}{\dhinterseq^{\gamma+1} w_n}=
\frac{1}{\gamma+1}\canditheta\;,
\end{align*}
with the corresponding statements for $t(N_1)$.
\end{corollary}

\subsubsection{Joint moments of the first and the last jump time}
To study joint moments we need a slightly stronger anticlustering condition \hyperlink{SummabilityAC}{$\mathcal{S}^{(0)}(\dhinterseq, \tepseq)$}.
Also, a function $f:\Rset_+^2\to \Rset_+$ is $(\gamma_1,\gamma_2)$-homogenous (with $\gamma_1,\gamma_2\geqslant 0$) if for all $a,b,s,t> 0$, $f(as,bt)=a^{\gamma_1}b^{\gamma_2}f(s,t)$. 
\begin{lemma}\label{lem:moments-first-and-last-jump}
Assume that \hyperlink{SummabilityAC}{$\mathcal{S}^{(0)}(\dhinterseq, \tepseq)$} holds. 
Let $f:\Rset_+^2\to \Rset_+$ be $(\gamma_1,\gamma_2)$-homogenous, $\gamma_1,\gamma_2\geqslant 0$. Then
\begin{align}\label{eq:joint-moments-locations}
\lim_{n\to\infty}\frac{1}{\dhinterseq^{\gamma_1+\gamma_2+1} w_n}\esp[f\left(t(1),t(N_1)\right)\ind{A_1}]=
\frac{f(1,1)}{\gamma_1+\gamma_2+1}\canditheta \;.
\end{align}
\end{lemma}

\subsection{Functionals dominated by the cluster length}\label{sec:extentions-cluster-length}
Note first that the cluster length functional defined in \eqref{eq:cluster-length-def} is unbounded, but shift-invariant. 
Recall the definition \eqref{eq:def-clusterlength} of the cluster length: $\clusterlength_{j}=t_{j}{(N_{j})}-t_{j}{(1)}+1$. 

We continue with an extension of \Cref{theo:cluster-RV}. 
We start with the following obvious result. We replace boundedness with the uniform integrability condition.  
\begin{proposition}
\label{theo:cluster-RV-unbounded}
  Assume that~\ref{eq:conditiondh} holds. Let $H$ be a continuous shift-invariant functional with support separated from $\bszero$. Assume that 
  \begin{align}\label{eq:vague-convergence-uniform-integrability}
  \lim_{\ell\to\infty}\limsup_{n\to\infty}\frac{\esp[ |H|(\tepseq^{-1}\bsX_{1,\dhinterseq})\ind{|H|(\tepseq^{-1}\bsX_{1,\dhinterseq})>\ell}]} {\dhinterseq w_n}=0\;. 
  \end{align}
  Then \eqref{eq:vague-convergence} holds. 
\end{proposition}
In this context, the main functional of interest will be a cluster length or any functional bounded by (a power of) the cluster length.

Recall \Cref{lem:moments-first-and-last-jump}. 
Since $\clusterlength_1-1=t(N_1)-t(1)$, one could expect that the limiting behaviour of the cluster length follows directly from that of the jump locations. However, we notice that the right of \eqref{eq:joint-moments-locations} vanishes for the function $f(x,y)=(y-x)_+$.  
  Indeed, the cluster length behaviour differs significantly. \Cref{lem:clusterlength-tight} (that follows directly from the vague convergence in \Cref{theo:cluster-RV}) gives tightness of the limiting conditional distribution of the (unscaled) cluster length. This is in contrast to the first and the last jump locations, that have to be scaled by the block size. 

We then extend \Cref{theo:cluster-RV} to the moments of the cluster length. Contrary to the jump times, there is a phase transition: for the first moment, the scaling is either $\dhinterseq w_n$ or $\dhinterseq^3 w_n^2$ depending on whether the blocks are small ($\dhinterseq^2 w_n\to 0$) or large ($\dhinterseq^2 w_n\to\infty$). In the case of the small blocks it is a consequence of vague convergence in \Cref{theo:cluster-RV} and the uniform integrability. See \Cref{lem:clusterlength-moments}. 
Here, a stronger anticlustering condition \ref{eq:conditionSstronger:gamma} is needed. We discuss its  optimality in \Cref{rem:optimality-of-rn2gamma}. \Cref{lem:clusterlength-moments-large} covers the case of large blocks.  As expected, a different scaling is needed.

\subsubsection{Distribution of the cluster length}
The following result is a re-statement of the one in \cite{drees:rootzen:2010} under weaker assumptions (see Lemma 5.2 and Eq.~(5.19) there). Only the anticlustering condition is needed here. The result follows directly from \Cref{theo:cluster-RV}.
\begin{lemma}\label{lem:clusterlength-tight}
Assume that \ref{eq:conditiondh} holds. Then for any finite $q\in\{0,\ldots,\dhinterseq\}$,
\begin{align}\label{eq:cluster-length-asymptotics}
\lim_{n\to\infty}
\frac{1}{\dhinterseq w_n}\pr\left(\{\clusterlength_1>q\}\cap A_1\right)
=\tailmeasurestar(\ind{\clusterlength>q})\;.
\end{align}
In particular,
the cluster length $\clusterlength_1$  is tight under the conditional law $\pr(\cdot\mid A_1)$.
\end{lemma}
We note that we can re-write \eqref{eq:cluster-length-asymptotics} in several equivalent ways:
\begin{align*}
\tailmeasurestar(\ind{\clusterlength>q})=\pr(\bsY_{-\infty,-1}^\ast\leqslant 1,\clusterlength(\bsY)>q)=
\pr(\bsY_{-\infty,-1}^\ast\leqslant 1,\bsY_{q,\infty}^\ast>1)=\canditheta \pr(\clusterlength(\bsZ)>q)\;.
\end{align*}
As will argue in \Cref{rem:clusterlength-properdistr}, the limiting distribution of the cluster length is a proper distribution. 
\subsubsection{Moments of the cluster length - small blocks scenario}
\Cref{lem:moments-first-and-last-jump} describes the limiting behaviour of the jump times. However, the
asymptotics for the moments of $(t(N_1)-t(1))$ changes completely. In order to proceed,
we need the stronger version of the anticlustering condition. Recall that \ref{eq:conditionSstronger:gamma} implies \ref{eq:conditiondh}. Also, recall that \ref{eq:conditionSstronger:gamma} forces the blocks to be small, since in the \iid\ case the anticlustering is equivalent to $\dhinterseq^{\gamma+1} w_n\to 0$.  As seen in \Cref{rem:optimality-of-rn2gamma} below, the latter small blocks condition is necessary for convergence of the moments of the cluster length.

The main result in case of small blocks is \Cref{lem:clusterlength-moments}. It follows from vague convergence along with the uniform integrability.  \Cref{lem:clusterlength-moments-large} covers the case of large blocks. 
\begin{lemma}\label{lem:cluster-length-varification}
Assume that 
\hyperlink{SummabilityAC}{$\mathcal{S}^{(\gamma)}(\dhinterseq, \tepseq)$} holds. 
Let $H\in\mch(\gamma)$. 
Then the uniform integrability condition \eqref{eq:vague-convergence-uniform-integrability} holds. 
\end{lemma}
We obtain immediately the following corollary.
\begin{corollary}\label{lem:clusterlength-moments}
Assume that 
\hyperlink{SummabilityAC}{$\mathcal{S}^{(\gamma+\delta)}(\dhinterseq, \tepseq)$} holds. 
Then for any $G\in\mch(\delta)$,
\begin{align*}
    %\label{eq:cluster-length-gamma-power-UniformIntegrability-interior}
   \lim_{n\to \infty} \frac{1}{\dhinterseq w_n} \esp \left[\clusterlength^\gamma(\tepseq^{-1}\bsX_{1,\dhinterseq}) G(\tepseq^{-1}\bsX_{1,\dhinterseq})\ind{ A_1 }\right]
     = \tailmeasurestar(G \clusterlength^\gamma)=\canditheta\esp \left[G(\bsZ)\clusterlength^\gamma(\bsZ) \right]\;.
\end{align*}
\end{corollary}
\subsubsection{Moments of the cluster length - large blocks scenario}
We recall that \hyperlink{SummabilityAC}{$\mathcal{S}^{(\gamma)}(\dhinterseq, \tepseq)$} is almost equivalent to $\dhinterseq^{\gamma+1}w_n\to 0$. Thus, we ask what happens if the latter condition is violated. Here, we need to impose some weak dependence assumptions on the underlying time series. Also, we state the result for $H=\clusterlength$. We were unable to establish result for a general functional $H$.
\begin{lemma}\label{lem:clusterlength-moments-large}
Assume that \ref{eq:conditiondh} holds and $\dhinterseq^{\gamma+1}w_n\to \infty$. 
Assume that $\bsX$ is
either $\ell$-dependent or $\beta$-mixing with the rates
\begin{align}\label{eq:mixing-rates-a}
\sum_{j=1}^{\dhinterseq} j^{\gamma}\beta_j=o(\dhinterseq^{\gamma+1}w_n^2)\;. 
\end{align}
Then 
\begin{align}\label{eq:moment-clusterlength-largeblocks}
&\lim_{n\to\infty}\frac{1}{\dhinterseq^{\gamma+2} w_n^2}\esp\left[(t(N_1)-t(1))^\gamma \ind{A_1}\right]=\lim_{n\to\infty}\frac{\esp\left[\clusterlength_1^\gamma \ind{A_1}\right]}{\dhinterseq^{\gamma+2} w_n^2}=\frac{1}{(\gamma+1)(\gamma+2)}\canditheta^2\;.
\end{align}
\end{lemma}
We accompany the above results with several remarks in \Cref{sec:additional-comments}.

\section{Weak convergence of empirical cluster process(es)}\label{sec:weak-convergence}
Let $H\in\mch(\gamma)$. 
The statistics 
\begin{align*}
\frac{1}{nw_n}\sum_{j=1}^{m_n}H(\mathbb{X}_j)=
\frac{1}{nw_n}\sum_{j=1}^{m_n}H(\tepseq^{-1}\bsX_{(j-1)\dhinterseq+1,j\dhinterseq})=:
\widetilde{\bsnu}_{n,\dhinterseq}^*(H)
\end{align*}
can be considered as a disjoint block pseudo-estimator of the cluster index $\tailmeasurestar(H)$. Indeed, under fairly general conditions one can establish 
\begin{align}\label{eq:conv-in-prob}
\widetilde{\bsnu}_{n,\dhinterseq}^*(H)\convprob \tailmeasurestar(H)\;. 
\end{align}
To obtain an estimator, set $k=nw_n$ and $\tepseq$ is replaced with the $k$th largest order statistics of $|\bsX_j|$, $j=1,\ldots,n$. Then, the consistency result still holds. See \cite[Section 10.1]{kulik:soulier:2020}. 

\subsection{Functionals dominated by the cluster length}
\subsubsection{Small blocks scenario}\label{sec:clt-small-blocks}
We assume $\dhinterseq^{\gamma+1}w_n\to 0$. Hence, \eqref{eq:conv-in-prob} is expected. We then consider the empirical cluster process defined in \eqref{eq:cluster-process}:
\begin{align*}
\widetilde{\mathbb{G}}_n(H)
=\sqrt{nw_n}\left\{\widetilde{\bsnu}_{n,\dhinterseq}^*(H)-\esp[\widetilde{\bsnu}_{n,\dhinterseq}^*(H)]\right\}
=\frac{1}{\sqrt{nw_n}}\sum_{j=1}^{m_n}\left\{H(\mathbb{X}_{j})-\esp[H(\mathbb{X}_{j})]\right\}\;,
\end{align*}
with the scaled blocks $\mathbb{X}_j$ defined in \eqref{eq:random-element-entire}.
To obtain its weak convergence, we need to impose a stronger condition on the block size. We will assume that \hyperlink{SummabilityAC}{$\mathcal{S}^{(\gamma(2+\delta))}(\dhinterseq, \tepseq)$} holds with an arbitrary small $\delta>0$, which in turn implies, in particular, that \eqref{eq:small-blocks-ep} holds. 

The following result extends Theorems 10.1.2 and 10.2.1 in \cite{kulik:soulier:2020} to unbounded functionals. We state this (and the next results) with all the assumptions, to be compatible with \cite{kulik:soulier:2020}.
\begin{theorem}\label{thm:weak-convergence}
Let $\{\bsX_j, j\in\Zset\}$ be a stationary, regularly varying $\Rset^d$-valued
time series. Assume that \ref{eq:rnbarFun0} and \ref{eq:conditionbeta-for-Aprime} hold. Fix $\gamma> 0$ and assume that there exists $\delta>0$ such that \hyperlink{SummabilityAC}{$\mathcal{S}^{(\gamma(2+\delta))}(\dhinterseq, \tepseq)$} holds. Then:
\begin{itemize}
\item For any $H\in \mch(\gamma)$, $\widetilde{\bsnu}_{n,\dhinterseq}^*(H)$ is a consistent pseudo-estimator of $\tailmeasurestar(H)$;
\item 
$\widetilde{\mathbb{G}}_n\convfidi \widetilde{\mathbb{G}}$ on $\mch(\gamma)$, where $\widetilde{\mathbb{G}}$ is a centered Gaussian process with the covariance functions $\tailmeasurestar(HH')$, $H,H'\in\mch(\gamma)$.    
\end{itemize}
\end{theorem}
We extend the above result to functional convergence. For this we need additional assumptions on the function class. We consider $\mcg\subseteq \mch(\gamma)$ and equip it with the metric $\rho^*(H,H')=\tailmeasurestar(\{H-H'\}^2)$. Define the random semi-metric $d_n$ on $\mcg$ by 
\begin{align*}
d_n^2(H,H')=\frac{1}{nw_n}\sum_{j=1}^{m_n}\{H\left(\mathbb{X}_j\right)-H'\left(\mathbb{X}_j\right)\}\;.
\end{align*}
In what follows, we need a random covering number of the function class. In the context of the stochastic processes considered here, we refer to Appendix C in \cite{kulik:soulier:2020}, while the general theory can be found in \cite{vandervaart:wellner:1996}. 
For statistical applications it suffices to consider linearly ordered function classes of the form $\mcg=\{H_s,s\in [s_0,t_0]\}$ for $0<s_0<t_0<\infty$ and some $H\in \mch(\gamma)$. For such classes  Conditions (b), (c), (d) below are automatically fulfilled. See e.g. the proof of Theorem 10.3.1 in \cite{kulik:soulier:2020}. 

The following result is a re-statement of Theorem C.4.5 in \cite{kulik:soulier:2020}, adapted to the setting of the present paper. See also Theorems 2.8 and 2.10 in \cite{drees:rootzen:2010}. A proof is not provided. 
\begin{theorem}\label{thm:weak-convergence-fclt}
Assume the conditions of \Cref{thm:weak-convergence} are satisfied. Assume also that: 
\begin{itemize}
\item[{\rm (a)}] The semi-metric space $(\mcg,\rho)$ is totally bounded. 
\item[{\rm (b)}] The envelope function $\bsG(\bsx):=\sup_{H\in \mcg}H(\bsx)$ belongs to $\mcg$.
\item[{\rm (c)}] For every sequence $\{\delta_n\}$ which decreases to zero, 
\begin{align*}
\lim_{n\to\infty}\sup_{H,H'\in \mcg; \rho(H,H')\leq \delta_n}\esp[d_n^2(H,H')]=0\;. 
\end{align*} 
\item[{\rm (d)}] There exists a measurable majorant of $N^*(\mcg,d_n, \epsilon)$ of the covering
number $N(\mcg, d_n, \epsilon)$ such that for every sequence $\{\delta_n\}$ which decreases
to zero,
\begin{align*}
\int_0^{\delta_n}\sqrt{\log N^*(\mcg,d_n, \epsilon)}\rmd \epsilon\convprob 0\;. 
\end{align*}
The the process $\{\mathbb{G}_n(H),H\in \mcg\}$ is asymptotically $\rho$-equicontinuous. \end{itemize}
\end{theorem}

\subsubsection{Moderate blocks scenario}\label{sec:clt-moderate-blocks}
We assume here that \eqref{eq:moderate-blocks-ep} holds.  
Hence, the consistency in \eqref{eq:conv-in-prob} is still expected.

Let $\mathbb{X}_j^\dag$, $j=1,\ldots,n$, be independent copies of $\mathbb{X}_j$ and define 
\begin{align*}
\widetilde{\mathbb{G}_n}^\dag(H):=
\frac{1}{\sqrt{nw_n}}\sum_{j=1}^{m_n}\left\{H(\mathbb{X}_{j}^\dag)-\esp[H(\mathbb{X}_{j})]\right\}\;.
\end{align*}
Thanks to the mixing assumption, the limiting behaviour of both $\widetilde{\mathbb{G}_n}^\dag$  and $\widetilde{\mathbb{G}_n}$ coincide. Set $H=\clusterlength^\gamma$. Then  
\begin{align*}
\var\left(\widetilde{\mathbb{G}_n}^\dag(H)\right)=\frac{\var(H(\mathbb{X}_1))}{\dhinterseq w_n}
=\underbrace{\frac{\esp[H^2(\tepseq^{-1}\bsX_{1,\dhinterseq})]}{\dhinterseq^{2\gamma+2}w_n^2}}_{=:J_{n,1}}
\left(\dhinterseq^{2\gamma+1}w_n\right)-\left(\underbrace{\frac{\esp[H(\tepseq^{-1}\bsX_{1,\dhinterseq})]}{\dhinterseq w_n}}_{=:{J_{n,2}}}\right)^2\dhinterseq w_n\;.
\end{align*}
On account of \Cref{lem:clusterlength-moments-large} and \Cref{lem:clusterlength-moments}, $J_{n,1}$ and $J_{n,2}$ converge to finite limits. (To be precise, convergence of $J_{n,2}$ requires the anticlustering condition, which is almost equivalent to $\dhinterseq^{\gamma+1}w_n\to 0$). Thus, under the moderate blocks condition the variance of the of the empirical cluster process diverges to infinity.  

Hence, the empirical process has to be modified. Define 
\begin{align*}
\widetilde{\mathbb{K}_n}(H):=\sqrt{\frac{n}{\dhinterseq^{2\gamma+1}}}
\left\{\widetilde{\bsnu}_{n,\dhinterseq}^*(H)-\esp[\widetilde{\bsnu}_{n,\dhinterseq}^*(H)]\right\}=
\frac{1}{\sqrt{n \dhinterseq^{2\gamma+1}}w_n}\sum_{j=1}^{m_n}\left\{H(\mathbb{X}_{j})-\esp[H(\mathbb{X}_{j})]\right\}\;.
\end{align*}
Then, with the obvious notation for the independent blocks process, 
\begin{align*}
\var\left(\widetilde{\mathbb{K}_n}^\dag(H)\right)=\frac{\var(H(\mathbb{X}_1))}{\dhinterseq^{2\gamma+2}w_n^2}
=\underbrace{\frac{\esp[H^2(\tepseq^{-1}\bsX_{1,\dhinterseq})]}{\dhinterseq^{2\gamma+2}w_n^2}}_{=J_{n,1}}
-\left(\underbrace{\frac{\esp[H(\tepseq^{-1}\bsX_{1,\dhinterseq})]}{\dhinterseq w_n}}_{={J_{n,2}}}\right)^2\dhinterseq^{-2\gamma}\;.
\end{align*}
Since $\gamma>0$ and $\dhinterseq\to\infty$, $\var\left(\widetilde{\mathbb{K}_n}^\dag(H)\right)$ converges as $n\to\infty$. We make this analysis precise in the following limiting statement. 
\begin{theorem}\label{thm:weak-convergence-moderate}
Let $\{\bsX_j, j\in\Zset\}$ be a stationary, regularly varying $\Rset^d$-valued
time series. Fix $\gamma>0$.  Assume that: 
\begin{itemize}
\item[{\rm (i)}] 
\ref{eq:rnbarFun0} holds.
\item[{\rm (ii)}]  $\bsX$ is either $\ell$-dependent or mixing with the rates \ref{eq:conditionbeta-for-Aprime} 
and \begin{align*}
\sum_{j=1}^{\dhinterseq} j^{(2+\delta)\gamma}\beta_j=o(\dhinterseq^{(2+\delta)\gamma+1}w_n^2)\;. 
\end{align*}
for some $\delta>0$. 
\item[{\rm (iii)}]  \hyperlink{SummabilityAC}{$\mathcal{S}^{(\gamma)}(\dhinterseq, \tepseq)$} holds (hence $\dhinterseq^{\gamma+1}w_n\to 0$), while $\dhinterseq^{2\gamma+1}w_n\to\infty$. 
\item[{\rm (iv)}]  
$n\dhinterseq w_n^2\to\infty$ and $\dhinterseq^{2\gamma+1}/n\to 0$ are satisfied.  
\end{itemize}
Then:
\begin{itemize}
\item For any $H\in \mch(\gamma)$, $\widetilde{\bsnu}_{n,\dhinterseq}^*(H)$ is a consistent pseudo-estimator of $\tailmeasurestar(H)$;
\item 
$\widetilde{\mathbb{K}}_n(\clusterlength^\gamma)\convdistr \widetilde{\mathbb{K}}(\clusterlength^\gamma)$, where $\widetilde{\mathbb{K}}(\clusterlength^\gamma)$ is a centered normal random variable with the variance 
$\{(2\gamma+1)(2\gamma+2)\}^{-1}\canditheta^2$. 
\end{itemize}
\end{theorem}
\subsubsection{Large blocks scenario}\label{sec:clt-large-blocks}
We finally assume that \eqref{eq:large-blocks-ep} holds. 
In view of \Cref{lem:clusterlength-moments}, the minimal condition required for consistency is that \hyperlink{SummabilityAC}{$\mathcal{S}^{(\gamma)}(\dhinterseq, \tepseq)$} is satisfied, which in turn yields $\dhinterseq^{\gamma+1}w_n\to 0$. If $\dhinterseq^{\gamma+1}w_n\to \infty$, then in view of \Cref{lem:clusterlength-moments-large} one can only expect that  
\begin{align*}
\frac{1}{n\dhinterseq^{\gamma+1}w_n^2}\sum_{j=1}^{m_n}\clusterlength^\gamma(\tepseq^{-1}\bsX_{(j-1)\dhinterseq+1,j\dhinterseq})
\end{align*}
converges in probability to the right hand side of \eqref{eq:moment-clusterlength-largeblocks}, which is not equal to $\tailmeasurestar(\clusterlength^\gamma)$. 
In particular, $\widetilde{\bsnu}_{n,\dhinterseq}^*(H)\convprob\infty$. 

\subsection{Location of jumps}\label{sec:jump-locations}
Let $G$ be bounded as in \Cref{lem:fist-jump-function-of}. 
Let $H=T_{\rm max}^\gamma G$, so that $H(\tepseq^{-1}\bsX_{1,\dhinterseq})=t(N_1)G(\tepseq^{-1}\bsX_{1,\dhinterseq})$, whenever there is a large jump in the block. In view of \Cref{lem:fist-jump-function-of}, the consistency result \eqref{eq:conv-in-prob} does not hold. The empirical cluster process has to be re-defined again. 
Let 
\begin{align*}
\widetilde{\widetilde{\bsnu}}_{n,\dhinterseq}^*(H):=
\frac{1}{n \dhinterseq^\gamma w_n}\sum_{j=1}^{m_n}H(\tepseq^{-1}\bsX_{(j-1)\dhinterseq+1,j\dhinterseq})\;.
\end{align*}
Thanks to \Cref{lem:fist-jump-function-of}, $\widetilde{\widetilde{\bsnu}}_{n,\dhinterseq}^*(H)$ may serve as a consistent estimator of $(\gamma+1)^{-1}\tailmeasurestar(G)$. 
In this spirit, we define another empirical cluster process:
\begin{align*}
\widetilde{\mathbb{L}}_n(H)
=\sqrt{nw_n}\left\{\widetilde{\widetilde{\bsnu}}_{n,\dhinterseq}^*(H)-
\esp\left[\widetilde{\widetilde{\bsnu}}_{n,\dhinterseq}^*(H)\right]\right\}
=\frac{1}{\sqrt{nw_n}\dhinterseq^\gamma}\sum_{j=1}^{m_n}\left\{H(\mathbb{X}_{j})-\esp[H(\mathbb{X}_{j})]\right\}\;.
\end{align*}
In other words, $\widetilde{\mathbb{L}}_n(H)=\widetilde{\mathbb{G}}_n(H_n)$, where $H_n=H/\dhinterseq^\gamma$. 

Then, with the obvious notation, 
\begin{align*}
\var\left(\widetilde{\mathbb{L}_n}^\dag(H)\right)=m_n\frac{\var(H(\mathbb{X}_1))}{n\dhinterseq^{2\gamma}
w_n}
=\underbrace{\frac{\esp[H^2(\tepseq^{-1}\bsX_{1,\dhinterseq})]}{\dhinterseq^{2\gamma+1}w_n}}_{=L_{n,1}}
-\left(\underbrace{\frac{\esp[H(\tepseq^{-1}\bsX_{1,\dhinterseq})]}{\dhinterseq^{\gamma+1} w_n}}_{={L_{n,2}}}\right)^2\dhinterseq w_n\;.
\end{align*}
By \Cref{lem:fist-jump-function-of}, the term $L_{n,1}$ converges to $(2\gamma+1)^{-1}\tailmeasurestar(G^2)$. Also, $L_{n,2}$ converges and hence the last term in the variance expansion vanishes. 
This motivates the following theorem. 
\begin{theorem}\label{thm:weak-convergence-jumps}
Let $\{\bsX_j, j\in\Zset\}$ be a stationary, regularly varying $\Rset^d$-valued
time series. Assume that \ref{eq:conditiondh}, \ref{eq:rnbarFun0} and \ref{eq:conditionbeta-for-Aprime} hold. Then:
\begin{itemize}
\item For any $G\in \mch(0)$, $\widetilde{\widetilde{\bsnu}}_{n,\dhinterseq}^*(T_{\rm max}^\gamma G)$ is a consistent pseudo-estimator of $(1+\gamma)^{-1}\tailmeasurestar(G)$;
\item 
$\widetilde{\mathbb{L}}_n\convfidi \widetilde{\mathbb{L}}$ on $\mch(0)$, where $\widetilde{\mathbb{L}}$ is a centered Gaussian process with the covariance function $(1+\gamma)^{-1}\tailmeasurestar(GG')$, $G,G'\in\mch(0)$.    
\end{itemize}
\end{theorem}
\begin{example}\label{xmpl:extremal-index-estimator}{\rm
Let $G(\bsx)=\ind{\bsx^*>1}$. In view of \eqref{eq:canditheta} below,  
\begin{align*}\widetilde{\theta}_n:=(\gamma+1)\widetilde{\widetilde{\bsnu}}_{n,\dhinterseq}^*(T_{\rm max}^\gamma G)\wedge 1
\end{align*} is a consistent pseudo-estimator the candidate extremal index $\canditheta$. This type of estimator can help with a finite sample downward bias - the classical blocks estimators are typically heavily biased towards zero. At the same time, $\widetilde{\theta}_n$ has typically a higher variance as compared to the classical blocks estimators. 
}
\end{example}

\section{Proofs I - unbounded cluster functionals}\label{sec:technical-details}
\subsection{Consequences of the anticlustering conditions}
\subsubsection{Formulas for the candidate extremal index}\label{sec:extremal-index}
Recall the notation $\bsx^*=\sup_{j\in\Zset}|\bsx_j|$. 
If $H(\bsx)=\ind{\bsx^\ast>1}$, since $|\bsY_0|>1$, \eqref{eq:cluster-measure} gives
\begin{align}\label{eq:canditheta}
\tailmeasurestar(\ind{\bsx^\ast>1})= \pr({\bsY_{-\infty,-1}^\ast\leqslant 1},{\bsY^\ast> 1})
=\pr({\bsY_{-\infty,-1}^\ast\leqslant 1})=\pr({\bsY_{1,\infty}^\ast\leqslant 1})=\canditheta\;.
\end{align}
\subsubsection{First and last jump decompositions}\label{sec:first-and-last-jump-decomposition}
The following decompositions will play a crucial role:
\begin{align*}
\ind{\bsX_{1,\dhinterseq}^\ast>\tepseq}=\sum_{i=1}^{\dhinterseq}
\ind{\bsX_{1,i-1}^\ast\leqslant \tepseq, \norm{\bsX_{i}}>\tepseq}\;.
\end{align*}
Recall the notation \eqref{eq:exceedance} for the exceedance times.
Set $t(i)=t_1{(i)}$, $i=1,\ldots,N_1$.
Then
\begin{align}\label{eq:decomposition-first-jump}
\{t(1)=j\}=\{\bsX_{1,j-1}^\ast\leqslant \tepseq, \norm{\bsX_{j}}>\tepseq\}\;,
\end{align}
(with the convention $\bsX_{1,0}^\ast\equiv 0$).
Likewise,
\begin{align}\label{eq:decomposition-last-jump}
\{t(N_1)=j\}=\{\bsX_{j+1,\dhinterseq}^\ast\leqslant \tepseq, \norm{\bsX_{j}}>\tepseq\}\;.
\end{align}
(with the convention $\bsX_{\dhinterseq+1,\dhinterseq}^\ast\equiv 0$). The above identities will be referred to as \textit{the first jump} and \textit{the last jump} decomposition, respectively.
\subsubsection{Conditional convergence of clusters}\label{sec:conditional-convergence-of-clusters}
As a consequence of \Cref{lem:tailprocesstozero}, for 
any bounded functional $H$ on $\ell_0(\Rset^d)$ that is almost
surely continuous with respect to the distribution of the tail process $\bsY$
any integers $i_1,i_2,j_1,j_2$ and arbitrary $s\geqslant 0$,
\begin{align}
\lim_{n\to\infty}\esp\left[H(\tepseq^{-1}(\bsX_{-[\dhinterseq s]-j_1,-i_1},\bsX_{i_2,j_2+[\dhinterseq s]})) \mid \norm{\bsX_0}>\tepseq\right]&=
\esp\left[H(\bsY_{-\infty,-i_1},\bsY_{i_2,\infty})\right]\;.\label{eq:tool-conditional-conv-rn-twosided}
\end{align}
This type of convergence will be referred to as the \textit{conditional convergence of clusters}.
\subsubsection{Vague convergence of clusters}\label{sec:integral-convergence-of-clusters}
In view of the anticlustering condition \ref{eq:conditiondh}, the following argument will be used repeatedly (cf. \Cref{theo:cluster-RV} and the proof of Theorem 6.2.5 in \cite{kulik:soulier:2020}). Let $G$ be a bounded cluster functional such that $G(\bsx)=0$ whenever $\bsx^\ast\leqslant 1$. Then
\begin{align}\label{eq:convergence-series-anticlustering-1}
\frac{1}{\dhinterseq}\sum_{i=1}^{\dhinterseq}\esp\left[G(\tepseq^{-1}\bsX_{1-i,\dhinterseq-i}) \ind{\bsX_{1-i,-1}^\ast\leqslant \tepseq} \mid\norm{\bsX_0}>\tepseq\right]
=\int_{0}^1 g_n(s)\rmd s
\end{align}
with
$$
g_n(s)=\esp\left[G(\tepseq^{-1}\bsX_{1-[\dhinterseq s],\dhinterseq-[\dhinterseq s]}) \ind{\bsX_{1-[\dhinterseq s],-1}^\ast\leqslant \tepseq} \mid\norm{\bsX_0}>\tepseq\right]
$$
converging to $g(s):= \esp[G(\bsY)\ind{\bsY_{-\infty,-1}^\ast\leqslant 1}]$ for each $s\in (0,1)$. By the dominated convergence the limit of the expression in \eqref{eq:convergence-series-anticlustering-1} becomes
\begin{align*}%\label{eq:convergence-series-anticlustering}
\esp[G(\bsY)\ind{\bsY_{-\infty,-1}^\ast\leqslant 1}]=\tailmeasurestar(G)\;.
\end{align*}
The first application is related to the candidate extremal index, cf. \eqref{eq:extremal-index}. 
Next, we extend \eqref{eq:convergence-series-anticlustering-1} in several directions. For $s\in (0,1)$ and integers $i_1,i_2\geqslant 1$, $0\leqslant j<i_2$, define
\begin{align*}%\label{eq:gn-twosided}
g_n^{\rm two-sided}(s;i_1,j,i_2):=\pr(\bsX_{1-[\dhinterseq s],-i_1}^\ast\leqslant \tepseq, \norm{\bsX_j}>\tepseq,\bsX_{i_2,\dhinterseq-[\dhinterseq s]}^*\leqslant \tepseq\mid \norm{\bsX_0}>\tepseq)\;.
\end{align*}
Then, \eqref{eq:tool-conditional-conv-rn-twosided} gives
\begin{align}\label{eq:gn-twosided-limit}
\lim_{n\to\infty} g_n^{\rm two-sided}(s;i_1,j,i_2)=\pr(\bsY_{-\infty,-i_1}^\ast\leqslant 1,\norm{\bsY_{j}}>1,\bsY_{i_2,\infty}^\ast\leqslant 1)\;
\end{align}
for all $s\in (0,1)$.
Likewise, let $G$ be a bounded cluster functional such that $G(\bsx)=0$ whenever $\bsx^\ast\leqslant 1$. Define
\begin{subequations}
\begin{alignat}{3}
g_n^{\rm left}(s;i_1,G)&:=\esp\left[G(\tepseq^{-1}\bsX_{1-[\dhinterseq s],\dhinterseq-[\dhinterseq s]})\ind{\bsX_{1-[\dhinterseq s],-i_1}^\ast\leqslant \tepseq} \mid \norm{\bsX_0}>\tepseq\right]\;, \label{eq:gn-onesided-left}\\
g_n^{\rm right}(s;i_2,G)&:= \esp\left[G(\tepseq^{-1}\bsX_{1-[\dhinterseq s],\dhinterseq-[\dhinterseq s]})\ind{\bsX_{i_2,\dhinterseq-[\dhinterseq s]}^\ast\leqslant \tepseq} \mid \norm{\bsX_0}>\tepseq\right]\;. \label{eq:gn-onesided-right}
\end{alignat}
\end{subequations}
Then for all $s\in (0,1)$,
\begin{subequations}
\begin{alignat*}{3}
\lim_{n\to\infty} g_n^{\rm left}(s;i_1,G)&=\esp[G(\bsY)\ind{\bsY_{-\infty,-i_1}^\ast\leqslant 1}]\;,
%\label{eq:gn-onesided-limit-left-general}
\\
\lim_{n\to\infty} g_n^{\rm right}(s;i_2,G)&=\esp[G(\bsY)\ind{\bsY_{i_2,\infty}^\ast\leqslant 1}]\;.%\label{eq:gn-onesided-limit-right-general}
\end{alignat*}
\end{subequations}
In particular, if $i_1=i_2=1$, then
\begin{align}\label{eq:gn-onesided-limit-genera-i1=21=1}
\lim_{n\to\infty} g_n^{\rm left}(s;1,G)=\lim_{n\to\infty} g_n^{\rm right}(s;1,G)=\tailmeasurestar(G)\;.
\end{align}
If moreover $G\equiv 1$, then we will simplify the notation $g_n^{\rm left}(s;1,G)$ as $g_n^{\rm left}(s)$ and
\begin{align}\label{eq:gn-onesided-limit}
\lim_{n\to\infty} g_n^{\rm left}(s)=\pr(\bsY_{-\infty,-1}^\ast\leqslant 1)=\canditheta\;, \ \ \lim_{n\to\infty} g_n^{\rm right}(s)=\pr(\bsY_{1,\infty}^\ast\leqslant 1)=\canditheta\;.
\end{align}
A function $f:\Rset_+\to \Rset_+$ is $\gamma$-homogeneous if for all $a,s>0$, $f(as)=a^{\gamma}f(s)$. Denote
\begin{align*}%\label{eq:f-onedim-scaled}
f_{n}(s)=f\left(\frac{[\dhinterseq s]}{\dhinterseq}\right)\;.
\end{align*}
Since a homogeneous function is continuous, we have
$
\lim_{n\to\infty}f_{n}(s)=f(s)$. 
Then, using \eqref{eq:gn-twosided-limit} and the bounded convergence,
\begin{align}
&\lim_{n\to\infty}\frac{1}{\dhinterseq^{\gamma+1}}
\sum_{i=1}^{\dhinterseq}f(i)\pr(\bsX_{1-i,-j_1}^\ast\leqslant \tepseq , \bsX_{j_2,\dhinterseq-i}^\ast\leqslant\tepseq\mid \norm{\bsX_0}>\tepseq)\notag\\
&=\lim_{n\to\infty}\int_0^1 f_n(s) g_n^{\rm two-sided}(s;j_1,0,j_2)\rmd s=\pr(\bsY_{-\infty,-j_1}^\ast\leqslant 1,\bsY_{j_2,\infty}^\ast\leqslant 1) \int_0^1 f(s)\rmd s\notag\\
&=\frac{f(1)}{\gamma+1}\pr(\bsY_{-\infty,-j_1}^\ast\leqslant 1,\bsY_{j_2,\infty}^\ast\leqslant 1)\;,\label{eq:integral-convergence-twosided}
\end{align}
Let now $G$ be a bounded cluster functional such that $G(\bsx)=0$ whenever $\bsx^\ast\leqslant 1$. Then, using \eqref{eq:gn-onesided-limit-genera-i1=21=1},
\begin{subequations}
\begin{alignat}{4}
&\lim_{n\to\infty}\frac{1}{\dhinterseq^{\gamma+1}}
\sum_{i=1}^{\dhinterseq}f(i)\esp\left[G(\bsX_{1-i,\dhinterseq-i})\ind{\bsX_{1-i,-1}^\ast\leqslant \tepseq \mid \norm{\bsX_0}>\tepseq}\right]\notag\\
&=\lim_{n\to\infty}\int_0^1 f_n(s) g_n^{\rm left}(s;1,G)\rmd s=\tailmeasurestar(G) \int_0^1 f(s)\rmd s=\frac{f(1)}{\gamma+1}\tailmeasurestar(G)\;,\label{eq:integral-convergence-left}\\
&\lim_{n\to\infty}\frac{1}{\dhinterseq^{\gamma+1}}
\sum_{i=1}^{\dhinterseq}f(i)\esp\left[G(\bsX_{1-i,\dhinterseq-i})\ind{\bsX_{1,\dhinterseq-i}^\ast\leqslant \tepseq \mid \norm{\bsX_0}>\tepseq}\right]\notag\\
&=\lim_{n\to\infty}\int_0^1 f_n(s) g_n^{\rm right}(s;1,G)\rmd s=\tailmeasurestar(G) \int_0^1 f(s)\rmd s=\frac{f(1)}{\gamma+1}\tailmeasurestar(G)\;.\label{eq:integral-convergence-right}
\end{alignat}
\end{subequations}
Next, we will say that a function $f:\Rset_+^2\to \Rset_+$ is a $(\gamma_1,\gamma_2)$-homogenous (with $\gamma_1,\gamma_2\geqslant 0$) if for all $a,b,s,t> 0$, $f(as,bt)=a^{\gamma_1}b^{\gamma_2}f(s,t)$. We will use the notation
\begin{align*}%\label{eq:f-scaled}
f_{n}(s,t)=f\left(\frac{[\dhinterseq s]}{\dhinterseq},\frac{[\dhinterseq t]}{\dhinterseq}\right)\;.
\end{align*}
Note that 
$\lim_{n\to\infty}f_{n}(s,t)=f(s,t)$. 
As in \eqref{eq:integral-convergence-twosided},
\begin{align}
&\lim_{n\to\infty}\frac{1}{\dhinterseq^{\gamma_1+\gamma_2+1}}
\sum_{i=1}^{\dhinterseq}f(i,i)\pr(\bsX_{1-i,-j_1}^\ast\leqslant \tepseq , \bsX_{j_2,\dhinterseq-i}^\ast\leqslant\tepseq\mid \norm{\bsX_0}>\tepseq)\notag\\
&=\lim_{n\to\infty}\int_0^1 f_n(s,s) g_n^{\rm two-sided}(s;j_1,0,j_2)\rmd s=\pr(\bsY_{-\infty,-j_1}^\ast\leqslant 1,\bsY_{j_2,\infty}^\ast\leqslant 1) \int_0^1 f(s,s)\rmd s\notag\\
&=\frac{f(1,1)}{\gamma_1+\gamma_2+1}\pr(\bsY_{-\infty,-j_1}^\ast\leqslant 1,\bsY_{j_2,\infty}^\ast\leqslant 1)\;.\label{eq:integral-multivariate-convergence-twosided}
\end{align}
Using \eqref{eq:gn-onesided-limit}, a similar argument gives
\begin{align}\label{eq:integral-multivariate-convergence-onesided}
&\lim_{n\to\infty}\frac{1}{\dhinterseq^{\gamma_1+\gamma_2+2}}
\sum_{i_1=1}^{\dhinterseq}\sum_{i_2=i_1+1}^{\dhinterseq}f(i_1,i_2)\pr(\bsX_{1-i_1,-1}^\ast\leqslant \tepseq\mid \norm{\bsX_{0}}>\tepseq)\pr(\bsX_{1,\dhinterseq-i_2}^\ast\leqslant \tepseq\mid \norm{\bsX_0}>\tepseq)\notag\\
&=
\lim_{n\to\infty}\int_0^1\left(\int_s^1 f_n(s,t) g_n^{\rm right}(t) \rmd t\right) g_n^{\rm left}(s)\rmd s
=\canditheta^2 f(1,1) \int_0^1s^{\gamma_1}\left(\int_s^1 t^{\gamma_2}\rmd t\right)\rmd s\;.
\end{align}
We will refer to this type of convergence as the \textit{vague convergence of clusters}.

\subsubsection{Proofs for unbounded cluster functionals}\label{sec:proof-for-internal}
\begin{proof}[Proof of \Cref{lem:first-and-last-jump}]
The first jump decomposition in \eqref{eq:decomposition-first-jump}, stationarity, vague convergence of clusters and the formula for $\canditheta$ in \eqref{eq:canditheta} give
\begin{align*}
&\frac{1}{\dhinterseq w_n}\pr(\{t(1)>q\dhinterseq\}\cap A_1)
=\frac{1}{\dhinterseq w_n}\sum_{i=q\dhinterseq+1}^{\dhinterseq}\pr(t(1)=i,\bsX_{1,\dhinterseq}^\ast>\tepseq)\\
&
= \frac{1}{\dhinterseq w_n}\sum_{i=q\dhinterseq+1}^{\dhinterseq}\pr(\bsX_{1,i-1}^\ast\leqslant \tepseq,\norm{\bsX_i}>\tepseq)
=\frac{1}{\dhinterseq}\sum_{i=q\dhinterseq+1}^{\dhinterseq}\pr(\bsX_{1-i,-1}^\ast\leqslant \tepseq\mid \norm{\bsX_0}>\tepseq)\\
&\to (1-q) \pr(\bsY_{-\infty,-1}^\ast\leqslant 1)=(1-q)\canditheta\;.
\end{align*}
The proof for $t(N_1)$ is analogous. We use the last jump decomposition \eqref{eq:decomposition-last-jump} instead.
\end{proof}
\begin{proof}[Proof of \Cref{lem:fist-jump-function-of}]
We use the first jump decomposition:
\begin{align*}
&\esp[f(t(1))G(\tepseq^{-1}\bsX_{1,\dhinterseq})]=\esp[f(t(1))G(\tepseq^{-1}\bsX_{1,\dhinterseq})\ind{A_1}]\\
&=w_n\sum_{i=1}^{\dhinterseq}f(i)\esp\left[G(\tepseq^{-1}\bsX_{1-i,\dhinterseq-i})\ind{\bsX_{1-i,-1}^\ast\leqslant \tepseq} \mid \norm{\bsX_0}>\tepseq\right]\;.
\end{align*}
Application of \eqref{eq:integral-convergence-left} yields the result.
The second part of the lemma follows from the decomposition on the last jump and \eqref{eq:integral-convergence-right}.
\end{proof}
\begin{proof}[Proof of \Cref{lem:moments-first-and-last-jump}]
We use the first and the last jump decomposition simultaneously:
\begin{align*}
&\esp[f(t(1),t(N_1))\ind{A_1}]\\
&=\sum_{i_1=1}^{\dhinterseq}\sum_{i_2=i_1}^{\dhinterseq}f(i_1,i_2)\pr(\bsX_{1,i_1-1}^\ast\leqslant \tepseq,\norm{\bsX_{i_1}}>\tepseq,\norm{\bsX_{i_2}}>\tepseq,\bsX_{i_2+1,\dhinterseq}^\ast\leqslant \tepseq)\\
%&=\sum_{i_1=1}^{\dhinterseq}f(i_1,i_1)\pr(\bsX_{1,i_1-1}^\ast\leqslant \tepseq,\norm{\bsX_{i_1}}>\tepseq,
%\bsX_{i_1+1,\dhinterseq}^\ast\leqslant \tepseq)\\
%&\phantom{=}+
%\sum_{i_1=1}^{\dhinterseq}\sum_{i_2=i_1+1}^{\dhinterseq}f(i_1,i_2)\pr(\bsX_{1,i_1-1}^\ast\leqslant %\tepseq,\norm{\bsX_{i_1}}>\tepseq,\norm{\bsX_{i_2}}>\tepseq,\bsX_{i_2+1,\dhinterseq}^\ast\leqslant \tepseq)\\
&=w_n\sum_{i=1}^{\dhinterseq}f(i,i)\pr(\bsX_{1-i,-1}^\ast\leqslant \tepseq,
\bsX_{1,\dhinterseq-i}^\ast\leqslant \tepseq\mid \norm{\bsX_0}>\tepseq)\\
&\phantom{=}+
\sum_{i=1}^{\dhinterseq}\sum_{j=1}^{\ell}f(i,j+i)\pr(\bsX_{1-i,-1}^\ast\leqslant \tepseq,\norm{\bsX_0}>\tepseq,\norm{\bsX_{j}}>\tepseq,\bsX_{j+1,\dhinterseq-i}^\ast\leqslant \tepseq)\\
&\phantom{=}+
\sum_{i_1=1}^{\dhinterseq}\sum_{i_2=i_1+\ell+1}^{\dhinterseq}f(i_1,i_2)\pr(\bsX_{1,i_1-1}^\ast\leqslant \tepseq,\norm{\bsX_{i_1}}>\tepseq,\norm{\bsX_{i_2}}>\tepseq,\bsX_{i_2+1,\dhinterseq}^\ast\leqslant \tepseq)\\
&=:J_1(\dhinterseq)+J_2(\dhinterseq;\ell)+\widetilde{J}_2(\dhinterseq;\ell)\;.
\end{align*}
Using \eqref{eq:integral-multivariate-convergence-twosided} we have
\begin{align}\label{eq:proof-of-moments-first-and-last-jump-1}
&\lim_{n\to\infty}\frac{J_1(\dhinterseq)}{\dhinterseq^{\gamma_1+\gamma_2+1}w_n}=
\frac{f(1,1)}{\gamma_1+\gamma_2+1}\pr(\bsY_{-\infty,-1}^\ast\leqslant 1,\bsY_{1,\infty}^\ast\leqslant 1)\;.
\end{align}
Set
\begin{align*}
D(\dhinterseq;j)=\sum_{i=1}^{\dhinterseq}f(i,j+i)\pr(\bsX_{1-i,-1}^\ast\leqslant \tepseq,\norm{\bsX_0}>\tepseq,\norm{\bsX_{j}}>\tepseq,\bsX_{j+1,\dhinterseq-i}^\ast\leqslant \tepseq)\;.
\end{align*}
Then
\begin{align*}
J_2(\dhinterseq;\ell)=\sum_{j=1}^{\ell}D(\dhinterseq;j)\;, \ \ \widetilde{J}_2(\dhinterseq;\ell)=\sum_{j=\ell+1}^{\dhinterseq}D(\dhinterseq;j)\;.
\end{align*}
For each $j=1,\ldots,\ell$, using the convergence \eqref{eq:gn-twosided-limit} we obtain
\begin{align*}
&\lim_{n\to\infty}\frac{1}{\dhinterseq^{\gamma_1+\gamma_2+1} w_n}D(\dhinterseq;j)=\lim_{n\to\infty}\int_0^1 f_n\left(s,\frac{j}{\dhinterseq}+s\right)g_n^{\rm two-sided}(s;1;j;j+1)\rmd s \\
&=\frac{f(1,1)}{\gamma_1+\gamma_2+1}\pr(\bsY_{-\infty,-1}^\ast\leqslant 1,\norm{\bsY_{j}}>1,\bsY_{j+1,\infty}^\ast\leqslant 1)\;.
\end{align*}
Hence,
\begin{align}\label{eq:proof-of-moments-first-and-last-jump-2}
\lim_{n\to\infty}\frac{1}{\dhinterseq^{\gamma_1+\gamma_2+1} w_n}J_2(\dhinterseq;\ell)=\frac{f(1,1)}{\gamma_1+\gamma_2+1}\sum_{j=1}^{\ell}\pr(\bsY_{-\infty,-1}^\ast\leqslant 1,\norm{\bsY_{j}}>1,\bsY_{j+1,\infty}^\ast\leqslant 1)\;.
\end{align}
Next, 
\begin{align*}
&\widetilde{J}_2(\dhinterseq;\ell)\leqslant
\constant\sum_{i_1=1}^{\dhinterseq}\sum_{i_2=i_1+\ell+1}^{\dhinterseq}i_1^{\gamma_1}i_2^{\gamma_2}\pr(\norm{\bsX_0}>\tepseq,\norm{\bsX_{i_2-i_1}}>\tepseq)\\
&=O(1)\sum_{i=1}^{\dhinterseq}i^{\gamma_1}\sum_{j=\ell+1}^{\dhinterseq-i}(j+i)^{\gamma_2}\pr(\norm{\bsX_0}>\tepseq,\norm{\bsX_{j}}>\tepseq)\\
&=O(\dhinterseq^{\gamma_1+1}\dhinterseq^{\gamma_2})\sum_{j=\ell+1}^{\dhinterseq}\pr(\norm{\bsX_0}>\tepseq,\norm{\bsX_{j}}>\tepseq)
\end{align*}
and hence
\begin{align*}
\frac{\widetilde{J}_2(\dhinterseq;\ell)}{\dhinterseq^{\gamma_1+\gamma_2+1}w_n}=O(1)\frac{1}{w_n}\sum_{j=\ell+1}^{\dhinterseq}\pr(\norm{\bsX_0}>\tepseq,\norm{\bsX_{j}}>\tepseq)
\end{align*}
and the latter term vanishes by letting $n\to\infty$ and then $\ell\to\infty$ on account of the stronger anticlustering condition \hyperlink{SummabilityAC}{$\mathcal{S}^{(0)}(\dhinterseq, \tepseq)$}.

Now, we combine \eqref{eq:proof-of-moments-first-and-last-jump-1}-\eqref{eq:proof-of-moments-first-and-last-jump-2} to get (keep in mind $\norm{\bsY_0}>1$):'
\begin{align*}
\frac{f(1,1)}{\gamma_1+\gamma_2+1}\sum_{j=0}^{\ell}\pr(\bsY_{-\infty,-1}^\ast\leqslant 1,\norm{\bsY_{j}}>1,\bsY_{j+1,\infty}^\ast\leqslant 1)\;.
\end{align*}
Letting $\ell\to\infty$ gives
\begin{align*}
\sum_{j=0}^{\infty}\pr(\bsY_{-\infty,-1}^\ast\leqslant 1,\norm{\bsY_{j}}>1,\bsY_{j+1,\infty}^\ast\leqslant 1)
=\pr(\bsY_{-\infty,-1}^\ast\leqslant 1,\bsY_{0,\infty}^\ast>1)=
%\pr(\bsY_{-\infty,-1}^\ast\leqslant 1)=
\canditheta\;.
\end{align*}
\end{proof}
\begin{proof}[Proof of \Cref{lem:clusterlength-tight}]
The result follows directly from \Cref{theo:cluster-RV}. Indeed, $H(\bsx)=\ind{\clusterlength(\bsx)>q}$ is bounded, has support separated from $\bszero$ and is continuous with respect to the law of $\bsY$. 

Moreover, \ref{eq:conditiondh} gives $\lim_{j\to\infty}|\bsY_j|=0$ almost surely, hence
\begin{align*}
\lim_{q\to\infty}\lim_{n\to\infty}
\pr\left(\clusterlength_1>q\mid A_1\right)
=0\;,
\end{align*}
that is the cluster length is tight under the conditional law.
\end{proof}
\begin{proof}[Proof of \Cref{lem:cluster-length-varification}]
As in the proof of \Cref{lem:moments-first-and-last-jump} we use the decompositions on the first and last jumps simultaneously. Since $H\in \mch(\gamma)$, we bound $|H|(\tepseq^{-1}\bsX_{1,\dhinterseq})\leqslant C(t(N_1)-t(1))^\gamma$.  Thus, 
\begin{align}\label{eq:proof-of-cluster-length-small}
&\esp\left[|H|(\tepseq^{-1}\bsX_{1,\dhinterseq})\ind{|H|(\tepseq^{-1}\bsX_{1,\dhinterseq})>L}\ind{A_1}\right]\nonumber\\
&\leqslant\constant \sum_{{i_1=1}\atop {i_2=i_1+L}}^{\dhinterseq}\esp\left[
|H|(\tepseq^{-1}\bsX_{i_1,i_2})\ind{\bsX_{1,i_1-1}^\ast\leqslant \tepseq,\norm{\bsX_{i_1}}>\tepseq,\norm{\bsX_{i_2}}>\tepseq,\bsX_{i_2+1,\dhinterseq}^\ast\leqslant \tepseq}\right]\nonumber\\
&=\constant\sum_{{i_1=1}\atop {i_2=i_1+L}}^{\dhinterseq}\esp\left[
|H|(\tepseq^{-1}\bsX_{0,i_2-i_1})\ind{\bsX_{1-i_1,-1}^\ast\leqslant \tepseq,\norm{\bsX_0}>\tepseq,\norm{\bsX_{i_2-i_1}}>\tepseq,\bsX_{i_2-i_1+1,\dhinterseq-i_1}^\ast\leqslant \tepseq}\right]\nonumber\\
&=\constant\sum_{i=1}^{\dhinterseq}\sum_{j=L}^{\dhinterseq-i}
\esp\left[|H|(\tepseq^{-1}\bsX_{0,j})\ind{\bsX_{1-i,-1}^\ast\leqslant \tepseq,\norm{\bsX_0}>\tepseq,\norm{\bsX_{j}}>\tepseq,\bsX_{j+1,\dhinterseq-i}^\ast\leqslant \tepseq}\right]\\
&\leqslant C\sum_{i=1}^{\dhinterseq}\sum_{j=L}^{\dhinterseq-i}
j^\gamma\pr(\bsX_{1-i,-1}^\ast\leqslant \tepseq,\norm{\bsX_0}>\tepseq,\norm{\bsX_{j}}>\tepseq,\bsX_{j+1,\dhinterseq-i}^\ast\leqslant \tepseq)\nonumber\\
&\leqslant C \dhinterseq\sum_{j=L}^{\dhinterseq}j^{\gamma}\pr(\norm{\bsX_0}>\tepseq,\norm{\bsX_{j}}>\tepseq)\;.\nonumber
\end{align}
We conclude by applying \hyperlink{SummabilityAC}{$\mathcal{S}^{(\gamma)}(\dhinterseq, \tepseq)$}.
\end{proof}
\begin{proof}[Proof of \Cref{lem:clusterlength-moments-large}]
We continue with \eqref{eq:proof-of-cluster-length-small} by setting $L=1$. Here, $H=\clusterlength$. Hence, 
\begin{align*}
&\esp\left[H(\tepseq^{-1}\bsX_{1,\dhinterseq}) \ind{A_1}\right]\\
&=\sum_{j=1}^{\dhinterseq} j^\gamma\sum_{i=1}^{\dhinterseq-j}
\pr\left(\bsX_{1-i,-1}^\ast\leqslant \tepseq,\norm{\bsX_0}>\tepseq,\norm{\bsX_{j}}>\tepseq,\bsX_{j+1,\dhinterseq-i}^\ast\leqslant \tepseq\right)\\
&=:\sum_{j=1}^{\dhinterseq} j^\gamma D(\dhinterseq;j)
=\sum_{j=1}^{\ell} j^\gamma D(\dhinterseq;j)+\sum_{j=\ell+1}^{\dhinterseq} j^\gamma D(\dhinterseq;j)\\
&=:J(\dhinterseq;\ell)+\widetilde{J}(\dhinterseq;\ell)\;.
\end{align*}
We have 
\begin{align*}
&\lim_{n\to\infty}\frac{1}{\dhinterseq w_n}D(\dhinterseq;j)\\
&=\lim_{n\to\infty} \int_0^1 \pr(\bsX_{1-[\dhinterseq s],-1}^\ast\leqslant \tepseq,\norm{\bsX_j}>\tepseq,\bsX_{j+1,\dhinterseq-[\dhinterseq s]}^\ast\leqslant \tepseq\mid \norm{\bsX_0}>\tepseq)\rmd s\\
&= \pr(
\bsY_{-\infty,-1}^\ast\leqslant 1,\norm{\bsY_{j}}>1,\bsY_{j+1,\infty}^\ast\leqslant 1)\;.
\end{align*}
Hence,  $J(\dhinterseq;\ell)=O_P(\dhinterseq w_n)$. 

Next, in case of  $\ell$-dependence, 
\begin{align*}
&\widetilde{J}(\dhinterseq;\ell)=w_n^2\sum_{j=\ell+1}^{\dhinterseq}j^\gamma \sum_{i=1}^{\dhinterseq-j}
\pr\left(\bsX_{1-i,-1}^\ast\leqslant\tepseq\mid \norm{\bsX_0}>\tepseq\right)\pr\left(\bsX^\ast_{1,\dhinterseq-i-j}\leqslant \tepseq\mid \norm{\bsX_0}>\tepseq\right)\\
&=:\widetilde{J}^*(\dhinterseq;\ell)\;.
\end{align*}
The term $\widetilde{J}^*(\dhinterseq;\ell)$ converges at the rate $\dhinterseq^{\gamma+2}w_n^2$. Indeed, using the notation \eqref{eq:gn-onesided-left}-\eqref{eq:gn-onesided-right} and the limit in \eqref{eq:gn-onesided-limit}, bearing in mind $[\dhinterseq s]+[\dhinterseq t]\sim [\dhinterseq (s+t)]$, we have
\begin{align*}
&\lim_{n\to\infty}\frac{1}{\dhinterseq^{\gamma+2}w_n^2}\esp\left[H(\tepseq^{-1}\bsX_{1,\dhinterseq})\ind{A_1}\right]\\
&=\lim_{n\to\infty}\int_0^1 t^{\gamma} \left(\int_0^{1-t}  g_n^{\rm left}(s) g_n^{\rm right}(s+t)\rmd s\right)\rmd t=\frac{\canditheta^2}{(\gamma+1)(\gamma+2)}\;.
\end{align*}
In case of $\beta$-mixing, 
\begin{align*}
&\widetilde{J}(\dhinterseq;\ell)=\widetilde{J}^*(\dhinterseq;\ell)+J_{\beta}(\dhinterseq)\;, 
\end{align*}
where 
\begin{align*}
J_{\beta}(\dhinterseq)\leq \constant 
\sum_{j=\ell+1}^{\dhinterseq} j^\gamma\sum_{i=1}^{\dhinterseq-j}\beta_j=
\constant 
\sum_{j=\ell+1}^{\dhinterseq} j^\gamma (\dhinterseq-j)\beta_j\leq \constant\ \dhinterseq 
\sum_{j=\ell+1}^{\dhinterseq} j^\gamma \beta_j
\;. 
\end{align*}
By the assumption, $J_{\beta}(\dhinterseq)=o(\dhinterseq^{\gamma+2}w_n^2)$. 
The proof for the large blocks is concluded.
\end{proof}
\subsubsection{Additional comments for unbounded cluster functionals}\label{sec:additional-comments}
\begin{remark}{\rm
The result of \Cref{lem:clusterlength-moments-large} should be compared to \Cref{lem:moments-first-and-last-jump}. Take $\gamma=\gamma_1=\gamma_2=1$ and $f(s,t)=(t-s)_+$. Then  
\Cref{lem:moments-first-and-last-jump} gives $\esp[(t(N_1)-t(1))\ind{A_1}]/(r_n^3 w_n)\to 0$. 
\Cref{lem:clusterlength-moments-large} allows to recover non-zero limit when scaled by $r_n^3 w_n^2$: $t(1)$ and $t(N_1)$ become independent when normalized by a smaller value. Indeed,  
assume that $U$ and $V$ are independent, standard uniform random variables. Then $\esp[(U-V)_+^\gamma]=1/((\gamma+1)(\gamma+2))$, which can be recognized as the right hand side of 
\eqref{eq:moment-clusterlength-largeblocks}.  
}\end{remark}
\begin{remark}{\rm 
The mixing condition \eqref{eq:mixing-rates-a} is somehow restrictive. Assume that $\beta_j=O(j^{-\delta})$, $\delta>0$. Then \eqref{eq:mixing-rates-a} holds whenever $\dhinterseq^{\delta}w_n^2\to \infty$. Thus, $\delta$ has to be sufficiently large. In particular, the condition holds if the mixing coefficients decay at the exponential rate.  
}
\end{remark}

\begin{remark}\label{rem:optimality-of-rn2gamma}
{\rm
The proof of tightness in \Cref{lem:clusterlength-tight} requires the anticlustering \ref{eq:conditiondh} only. The convergence of the moments requires the stronger condition  \ref{eq:conditionSstronger:gamma}. In \iid\ case, the stronger anticlustering condition is equivalent to $\dhinterseq^{\gamma+1}w_n=o(1)$; cf. \Cref{rem:rn2gamma}. We argue that the latter assumption is necessary for convergence of moments of $\clusterlength_1$. 

Assume that $\bsX$ is \iid\  Then,
the cluster length of the tail process is $\clusterlength(\bsY) \stackrel{\text{d}}{=} 1$. We find sufficient and necessary conditions so that
\begin{align} \label{eq;lim-clusterlen1}
 \lim_{n\to \infty} \esp\left[ \clusterlength_1 \mid A_1 \right] = 1\;.
\end{align}
We have as $n\to\infty$,
\begin{align*}
 f(1) &:=  \pr\left(  \clusterlength_1 = 1 \mid A_1 \right) = \binom{\dhinterseq}{1} w_n (1-w_n)^{\dhinterseq - 1} \bigg/  \Big[   1 - (1-w_n)^{\dhinterseq}  \Big] \sim 1.
%\label{eq;clusterlen-pmf-q=1;1}
\end{align*}
Next, for any $i \geqslant 2$, we have
\begin{align*}
    f(i) &:=  \pr\left( \clusterlength_1 = i \mid A_1 \right) \notag \\
    &=  \frac{1}{  1 - (1-w_n)^{\dhinterseq} }
    \sum_{j=1} ^ {\dhinterseq - i + 1}
    \underbrace{ ( 1-w_n)^{j-1}  w_n }_{ t(1) = j }  \cdot
    \underbrace{( 1 - w_n ) ^ {\dhinterseq - (j+i-1)} w_n}_{ t(N_1) = j+i-1  } \notag \\
    &= \frac{ w_n^2 (\dhinterseq-i+1) (1-w_n)^{\dhinterseq-i} }{1 - (1-w_n)^{\dhinterseq}}
    \sim  \frac{\dhinterseq-i+1}{\dhinterseq} w_n.
    %\label{eq;clusterlen-pmf-q>1;1}
\end{align*}
uniformly for all $2\leqslant i \leqslant \dhinterseq$.

Clearly, the limit in \eqref{eq;lim-clusterlen1} is equivalent to $\sum_{i=2}^{\dhinterseq} i^\gamma f(i) \to 0$.
It suffices to note that \begin{align*}
    \sum_{i=2}^{\dhinterseq} i^\gamma f(i)
    &\sim
    \sum_{i=2}^{\dhinterseq} i^\gamma w_n  \frac{\dhinterseq - i+1}{\dhinterseq} \\
     &= \dhinterseq ^ {\gamma + 1} w_n \sum_{i=2}^{\dhinterseq} \left(
     \frac{i}{\dhinterseq} \right)^ \gamma \cdot \frac{\dhinterseq - i +1}{\dhinterseq} \cdot \frac{1}{\dhinterseq}
     \sim  \dhinterseq^{\gamma + 1} w_n {\rm Beta}(\gamma+1,2)\;,
\end{align*}
where ${\rm Beta}(\cdot, \cdot)$ is the standard Beta function. This yields that the condition $\dhinterseq^{\gamma+1}w_n\to 0$ is necessary (and sufficient) for the conditional convergence of moments of the cluster length.
}
\end{remark}
\begin{remark}\label{rem:clusterlength-properdistr}{\rm 
From \eqref{eq:cluster-length-asymptotics} we obtain for $q\geqslant 1$,
\begin{align*}%\label{eq:cluster-length-asymptotics-prob}
\lim_{n\to\infty}
\pr\left(\clusterlength_1=q\mid A_1\right)
=\frac{\pr(\bsY_{-\infty,-1}^\ast\leqslant 1,\norm{\bsY_{q-1}}>1,\bsY_{q,\infty}^\ast\leqslant 1)}{\canditheta}=:f(q)\;.
\end{align*}
Then, since $\norm{\bsY_0}>1$, bearing in mind \eqref{eq:canditheta},
\begin{align*}
f(1)+\sum_{q=2}^\infty f(q)&=
\frac{\pr(\bsY_{-\infty,-1}^\ast\leqslant 1,\bsY_{1,\infty}^\ast\leqslant 1)}{\canditheta}
+\frac{1}{\canditheta}\sum_{q=2}^\infty
{\pr(\bsY_{-\infty,-1}^\ast\leqslant 1,\norm{\bsY_{q-1}}>1,\bsY_{q,\infty}^\ast\leqslant 1)}\\
&=\frac{\pr(\bsY_{-\infty,-1}^\ast\leqslant 1,\bsY_{1,\infty}^\ast\leqslant 1)}{\canditheta}+\frac{\pr(\bsY_{-\infty,-1}^\ast\leqslant 1,\bsY_{1,\infty}^\ast> 1)}{\canditheta}=1\;.
\end{align*}
Hence, the limiting conditional distribution of the cluster length is a proper distribution.}
\end{remark}

\begin{remark}{\rm 
Take $q=0$. The the left hand side of \eqref{eq:cluster-length-asymptotics} is by the definition of the cluster length equal to 1. At the same time, since $|\bsY_0|>1$, using \eqref{eq:canditheta}, the right hand side of \eqref{eq:cluster-length-asymptotics} also reduces to one.
In case of extremal independence, $\bsY_{1,\infty}^\ast\equiv 0$ and hence $\lim_{n\to\infty}\pr(\clusterlength_1>1\mid A_1)=0$ as expected.
}
\end{remark}
\section{Proofs II - the empirical cluster processes}\label{sec:proofs-clt}
\begin{proof}[Proof of \Cref{thm:weak-convergence}]
We start with consistency. 
We note that \hyperlink{SummabilityAC}{$\mathcal{S}^{(\gamma(2+\delta))}(\dhinterseq, \tepseq)$} implies 
\hyperlink{SummabilityAC}{$\mathcal{S}^{(\gamma)}(\dhinterseq, \tepseq)$}. Hence 
\Cref{lem:clusterlength-moments} gives $\lim_{n\to\infty}\esp\left[\widetilde{\bsnu}_{n,\dhinterseq}^*(H)\right]=\tailmeasurestar(H)$. 
Thanks to the mixing assumption, we can treat the blocks as independent. Then 
\begin{align*}
\var\left(\widetilde{\bsnu}_{n,\dhinterseq}^*(H)\right)\leq m_n\frac{\esp[H^2(\tepseq^{-1}\bsX_{1,\dhinterseq})]}{n^2w_n^2}=\frac{\esp[H^2(\tepseq^{-1}\bsX_{1,\dhinterseq})]}{\dhinterseq w_n}\frac{1}{nw_n}\;. 
\end{align*}
Since \hyperlink{SummabilityAC}{$\mathcal{S}^{(2\gamma)}(\dhinterseq, \tepseq)$} holds, by \Cref{lem:clusterlength-moments} and $nw_n\to\infty$, the variance of the estimator vanishes. Consistency has been established. 

Next,
we verify that the class $\mch(\gamma)$ fulfills the conditions of Theorem 10.2.1 in \cite{kulik:soulier:2020}. 
\begin{itemize}
\item[{\rm (BCLT1)}] This conditions states that \eqref{eq:thelimitwhichisnolongercalledbH} holds. This in turn is guaranteed by \Cref{lem:cluster-length-varification}. 
\item[{\rm (BCLT2)}] We need to show that for any $\eta>0$ we have
\begin{align*}
\lim_{n\to\infty}\frac{1}{\dhinterseq w_n}\esp[H^2(\tepseq^{-1}\bsX_{1,\dhinterseq})\ind{|H|(\tepseq^{-1}\bsX_{1,\dhinterseq})>\eta \sqrt{nw_n}}]=0\;. 
\end{align*} 
The left hand-side is bounded by 
\begin{align*}
\lim_{n\to\infty}\frac{\esp[|H|^{2+\delta}(\tepseq^{-1}\bsX_{1,\dhinterseq})]}{\dhinterseq w_n}\frac{1}{(nw_n)^{\delta/2}}
\end{align*}
which vanishes on account of the anticlustering assumption and $nw_n\to\infty$. Here, the anticlustering condition with $\gamma(2+\delta)$ is needed. Indeed, if $H\in \mch(\gamma)$, then $|H|^{2+\delta}$ belongs to $\mch(\gamma(2+\delta))$. 
\item[{\rm (BCLT3)}]
We need to verify that there exist functions $K_n:(\Rset^d)^{\ell_n}\to \Rset$ such that 
\begin{align*}
\left|H(\tepseq^{-1}\bsX_{1,\dhinterseq})-H(\tepseq^{-1}\bsX_{1,\dhinterseq-\ell_n})\right|\leq K_n(\bsX_{\dhinterseq-\ell_n+1,\dhinterseq})
\end{align*}
and 
\begin{align*}
\lim_{n\to\infty}\frac{\esp[K_n(\bsX_{\dhinterseq-\ell_n+1,\dhinterseq})]}{\dhinterseq w_n}=0\;. 
\end{align*}
On account of
\Cref{Assumption:class-mathcalH}(iii), $H(\tepseq^{-1}\bsX_{1,\dhinterseq})\not=H(\tepseq^{-1}\bsX_{1,\dhinterseq-\ell_n})$ if and only if 
$\bsX^*_{\dhinterseq-\ell_n+1,\dhinterseq}>\tepseq$. Thus, choose $K_n(\bsX_{\dhinterseq-\ell_n+1,\dhinterseq})=\ind{\bsX^*_{\dhinterseq-\ell_n+1,\dhinterseq}>\tepseq}$ and then 
\begin{align*}
\lim_{n\to\infty}\frac{\pr\left(\bsX^*_{\dhinterseq-\ell_n+1,\dhinterseq}>\tepseq\right)}{\dhinterseq w_n}=
\lim_{n\to\infty}\frac{\pr\left(\bsX^*_{1,\ell_n}>\tepseq\right)}{\ell_n w_n}\frac{\ell_n w_n}{\dhinterseq w_n}=0\;. 
\end{align*}
\end{itemize}
\end{proof}
\begin{proof}[Proof of \Cref{thm:weak-convergence-moderate}]
The proof is almost identical to the one for \Cref{thm:weak-convergence}. 

Since 
\hyperlink{SummabilityAC}{$\mathcal{S}^{(\gamma)}(\dhinterseq, \tepseq)$} holds,
\Cref{lem:clusterlength-moments} gives $\lim_{n\to\infty}\esp\left[\widetilde{\bsnu}_{n,\dhinterseq}^*(H)\right]=\tailmeasurestar(H)$. 
Thanks to the mixing assumption, we can treat the blocks as independent. Then 
\begin{align*}
\var\left(\widetilde{\bsnu}_{n,\dhinterseq}^*(H)\right)\leq m_n\frac{\esp[H^2(\tepseq^{-1}\bsX_{1,\dhinterseq})]}{n^2w_n^2}=
\frac{\esp[H^2(\tepseq^{-1}\bsX_{1,\dhinterseq})]}{\dhinterseq^{2\gamma+2} w_n^2}\frac{\dhinterseq^{2\gamma+1}}{n}\;. 
\end{align*}
Since $\dhinterseq^{2\gamma+1} w_n\to\infty$, we use 
\Cref{lem:clusterlength-moments-large} to conclude convergence of the second last ratio. 
At the same time. 
$\dhinterseq^{2\gamma+1}/n\to 0$ by the assumption. Consistency has been established.

For the weak convergence, 
only the modified version of the Lindeberg condition (BCLT2) has to be proven. 
We need to show that for any $\eta>0$ we have
\begin{align*}
\lim_{n\to\infty}\frac{1}{\dhinterseq^{2\gamma+2} w_n^2}\esp\left[H^2(\tepseq^{-1}\bsX_{1,\dhinterseq})\ind{|H|(\tepseq^{-1}\bsX_{1,\dhinterseq})>\eta \sqrt{n\dhinterseq^{2\gamma+1}}w_n}\right]=0\;. 
\end{align*} 
The left hand-side is bounded by 
\begin{align*}
\lim_{n\to\infty}\frac{\esp[|H|^{2+\delta}(\tepseq^{-1}\bsX_{1,\dhinterseq})]}{\dhinterseq^{(2+\delta)\gamma+2} w_n^2}\frac{1}{(n\dhinterseq w_n^2)^{\delta/2}}\;.
\end{align*}
The first ratio converges by \Cref{lem:clusterlength-moments-large} (thanks to the assumptions (ii)-(iii)), while the second one 
vanishes by the assumption (iv). 
\end{proof} 
\begin{proof}[Proof of \Cref{thm:weak-convergence-jumps}]
Thanks to \Cref{lem:fist-jump-function-of}, 
\begin{align*}
\lim_{n\to\infty}
\esp\left[\widetilde{\widetilde{\bsnu}}_{n,\dhinterseq}^*(H)\right]=\lim_{n\to\infty}
\frac{\esp[T_{\rm max}^\gamma(\tepseq^{-1}\bsX_{1,\dhinterseq})G(\tepseq^{-1}\bsX_{1,\dhinterseq})]}{\dhinterseq^{\gamma+1}w_n}=\frac{1}{\gamma+1}\tailmeasurestar(G)\;. 
\end{align*}
Again, the blocks can be treated as independent and hence
\begin{align*}
\var\left(\widetilde{\widetilde{\bsnu}}_{n,\dhinterseq}^*(H)\right)\leq
\frac{\esp[T_{\rm max}^{2\gamma}(\tepseq^{-1}\bsX_{1,\dhinterseq})G^2(\tepseq^{-1}\bsX_{1,\dhinterseq})]}{\dhinterseq^{2\gamma+1}w_n}
\frac{1}{nw_n}\to 0\;.  
\end{align*}
For the weak convergence, again, 
only the modified version of the Lindeberg condition (BCLT2) has to be proven. 
We need to show that for any $\eta>0$ we have
\begin{align*}
\lim_{n\to\infty}\frac{1}{\dhinterseq^{2\gamma+1} w_n}\esp\left[H^2(\tepseq^{-1}\bsX_{1,\dhinterseq})\ind{|H|(\tepseq^{-1}\bsX_{1,\dhinterseq})>\eta \sqrt{n\dhinterseq^{2\gamma}w_n}}\right]=0\;. 
\end{align*} 
The left hand-side is bounded by 
\begin{align*}
\lim_{n\to\infty}\frac{\esp[|H|^{2+\delta}(\tepseq^{-1}\bsX_{1,\dhinterseq})]}{\dhinterseq^{(2+\delta)\gamma+1} w_n}\frac{1}{(nw_n)^{\delta/2}}\;
\end{align*}
and vanishes by \Cref{lem:fist-jump-function-of} and $nw_n\to\infty$. 
\end{proof}
\section*{Acknowledgements}
The authors are grateful to Philippe Soulier for an extensive discussion.  
Both authors were supported by an NSERC grant. A part of this research was conducted during the second author's stay at the Sydney Mathematical Research Institute (October - November 2022). Rafa{\l} Kulik thanks for the support and hospitality provided by SMRI.

%\section{Bibliography}
%\bibliographystyle{plain}


\begin{thebibliography}{10}

\bibitem{basrak:planinic:soulier:2018}
Bojan Basrak, Hrvoje Planini\'{c}, and Philippe Soulier.
\newblock An invariance principle for sums and record times of regularly
  varying stationary sequences.
\newblock {\em Probability Theory and Related Fields}, 172(3-4):869--914, 2018.

\bibitem{basrak:segers:2009}
Bojan Basrak and Johan Segers.
\newblock Regularly varying multivariate time series.
\newblock {\em Stochastic Processes and their Applications}, 119(4):1055--1080,
  2009.

\bibitem{bucher:zhou:2018}
Axel B\"{u}cher and Chen Zhou.
\newblock A horse racing between the block maxima method and the peak over
  threshold approach.
\newblock arXiv:1807.00282, 2018.

\bibitem{chen:kulik:2023b}
Zaoli Chen and Rafal Kulik.
\newblock Asymptotic expansions for blocks estimators: Pot framework.
\newblock In preparation.

\bibitem{davis:hsing:1995}
Richard~A. Davis and Tailen Hsing.
\newblock Point process and partial sum convergence for weakly dependent random
  variables with infinite variance.
\newblock {\em Annals of Probability}, 23(2):879--917, 1995.

\bibitem{drees:rootzen:2010}
Holger Drees and Holger Rootz{\'e}n.
\newblock Limit theorems for empirical processes of cluster functionals.
\newblock {\em Annals of Statistics}, 38(4):2145--2186, 2010.

\bibitem{hsing:1991}
Tailen Hsing.
\newblock Estimating the parameters of rare events.
\newblock {\em Stochastic Processes and their Applications}, 37(1):117--139,
  1991.

\bibitem{hsing:1993}
Tailen Hsing.
\newblock Extremal index estimation for a weakly dependent stationary sequence.
\newblock {\em Annals of Statistics}, 21(4):2043--2071, 1993.

\bibitem{kulik:soulier:2020}
Rafa{\l} Kulik and Philippe Soulier.
\newblock {\em Heavy tailed time series}.
\newblock Springer, 2020.

\bibitem{last:2023}
G\"{u}nter Last.
\newblock Tail processes and tail measures: An approach via palm calculus.
\newblock {\em Extremes}, 2023.

\bibitem{mikosch:wintenberger:2016}
Thomas Mikosch and Olivier Wintenberger.
\newblock A large deviations approach to limit theorem for heavy-tailed time
  series.
\newblock {\em Probability Theory and Related Fields}, 166(1-2):233--269, 2016.

\bibitem{planinic:2023}
Hrvoje Planini\'{c}.
\newblock Palm theory for extremes of stationary regularly varying time series
  and random fields.
\newblock {\em Extremes}, 26:45--82, 2023.

\bibitem{planinic:soulier:2018}
Hrvoje Planini\'{c} and Philippe Soulier.
\newblock The tail process revisited.
\newblock {\em Extremes}, 21(4):551--579, 2018.

\bibitem{robert:2009}
Christian~Y. Robert.
\newblock Inference for the limiting cluster size distribution of extreme
  values.
\newblock {\em Annals of Statistics}, 37(1):271--310, 2009.

\bibitem{rootzen:leadbetter:dehaan:1998}
Holger Rootz{\'e}n, Ross~M. Leadbetter, and Laurens de~Haan.
\newblock On the distribution of tail array sums for strongly mixing stationary
  sequences.
\newblock {\em Annals of Applied Probability}, 8(3):868--885, 1998.

\bibitem{smith:weissman:1994}
Richard~L. Smith and Ishay Weissman.
\newblock Estimating the extremal index.
\newblock {\em Journal of the Royal Statistical Society. Series B.
  Methodological}, 56(3):515--528, 1994.

\bibitem{vandervaart:wellner:1996}
Aad~W. van~der Vaart and Jon~A. Wellner.
\newblock {\em Weak convergence and empirical processes}.
\newblock Springer, New York, 1996.

\end{thebibliography}
\end{document}